\documentclass[11pt]{article}
\usepackage{amsfonts,color}
\usepackage{amsmath,amssymb,tipa,enumerate,txfonts,mathrsfs,algorithm,float,indentfirst}
\usepackage{bm}
\usepackage{a4wide}
\usepackage{epsfig,multirow,epstopdf}
\usepackage{graphicx}

\makeatletter
  \newcommand\figcaption{\def\@captype{figure}\caption}
  \newcommand\tabcaption{\def\@captype{table}\caption}
\makeatletter

\newtheorem{theorem}{Theorem}[section]
\newtheorem{remark}{Remark}[section]

\def\be{\begin{equation}}
\def\en{\end{equation}}
\def\beq{\begin{eqnarray}}
\def\eq{\end{eqnarray}}
\def\beqx{\begin{eqnarray*}}
\def\eqx{\end{eqnarray*}}

\def\12{{1\over 2}}

\def\0{{\bf 0}}

\def\prod{{\bf \sqcap}}

\def\Pi{\pi}

\begin{document}
\title{A framework of discontinuous Galerkin neural networks for iteratively approximating residuals}
\author{
Long Yuan\footnote{Corresponding author. College of Mathematics and Systems Science, Shandong University of Science and Technology,
Qingdao 266590, China. This author was supported by Shandong Provincial Natural Science Foundation under the grant ZR2024MA059. (email:\it{sdbjbjsd@163.com}).
%The author worked on the details of the algorithms and analysis, and carried out the numerical experiments.
}
~~~and~~~
Hongxing Rui\footnote{
Corresponding author. School of Mathematics, Shandong University, Jinan 250100, China (hxrui@sdu.edu.cn).
%The author proposed the ideas of the method and analysis.
}
}

\date{\today }
\maketitle

\begin{abstract} We propose an abstract discontinuous Galerkin neural network (DGNN) framework for analyzing the convergence of least-squares
methods based on the residual minimization when feasible solutions are neural networks. Within this framework, we define a quadratic loss functional as in the least square method with $h-$refinement and introduce new discretization sets spanned by element-wise neural network functions. The desired neural network approximate solution is recursively supplemented by solving a sequence of quasi-minimization problems associated with the underlying loss functionals and the adaptively augmented discontinuous neural network sets without the assumption on the boundedness of the neural network parameters. We further propose a discontinuous Galerkin Trefftz neural network discretization (DGTNN) only with a single hidden layer to reduce the computational costs. Moreover, we design a template based on the considered models for initializing nonlinear weights. Numerical experiments confirm that compared to existing PINN algorithms, the proposed DGNN method with one or two hidden layers is able to improve the relative $L^2$ error by at least one order of magnitude at low computational costs.
\end{abstract}


\vskip 0.1in


\section{Introduction}
The neural network method has emerged as a powerful tool for approximating solutions to partial differential equations (PDEs). To obtain an approximate solution of a PDE using neural networks, the key step is indeed to minimize the PDE residual effectively. Several approaches have been proposed to achieve this \cite{berg, eyu, JKK, hexu, Kharazmi, raissi, zang}. The work \cite{JKK} proposes a conservative physics-informed neural network (cPINN) on the discrete sub-domains for nonlinear conservation laws, where PINN is applied and the conservation property of cPINN is achieved by enforcing the numerical flux continuity in the {\it strong} form along the sub-domain interfaces. Besides, a general framework for $hp$-variational physics-informed neural networks ($hp$-VPINNs) is formulated in \cite{KZK}, where the trial space consisting of neural network sets is defined {\it globally} over the entire computational domain, while the test space contains {\it piecewise} polynomials. Specifically the $hp$-refinement in this study corresponds to a global approximation with a local learning algorithm. While such physics-informed neural networks and variational approaches have enjoyed success in particular cases, {\it the stagnation of relative error around ${\bf 0.1-1.0\%}$ still exists,} no matter how many neurons or layers are used to define the underlying network architecture. Pushing beyond this limit requires addressing optimization, architecture and sampling challenges.

%new path
While along another distinct school of thought, {\it continuous} neural networks combined with the Galerkin framework, based on the adaptive construction of a sequence of finite-dimensional subspaces whose basis functions are realizations of a sequence of deep neural networks, have been explored in \cite{AD}. The sequential nature of the algorithm provides a disciplined framework to achieving the accuracy of a given approximation.
%The sequential nature of the algorithm offers a systematic approach to enhancing the accuracy of a given approximation.
At the same time, this approach relies on the existing assumption that the bilinear form of the associated variational problem satisfies {\it symmetric, bounded} and {\it coercive} properties, and the convergence results established rely on a key assumption that the neural network parameters are bounded.

Recently, a {\it discontinuous} Galerkin plane wave neural network (DGPWNN) with $hp$-refinement and a single hidden layer for approximately solving the {\bf \it homegeneous} Helmholtz equation and time-harmonic Maxwell equations with {\bf \it constant coefficients} was proposed in \cite{ywh}, where the adaptive construction of recursively augmented discontinuous Galerkin subspaces is introduced, and their basis functions approaching the unit residuals are recursively generated by iteratively solving quasi-maximization problems associated with the underlying residual functionals and the intersection of the closed unit ball and discontinuous plane wave neural network spaces. Besides this work, the very recent work \cite{yuanhu} proposed a discontinuous plane wave neural network (DPWNN), where the desired approximate solution is recursively generated by iteratively solving the quasi-minimization problem associated with the quadratic functional and the sets spanned by element-wise plane wave neural network functions with a single hidden layer. Both works \cite{ywh,yuanhu} employ the plane wave functions as activation functions on each element, which seriously hinders the extension of the method to {\it nonhomogeneous} models with {\bf \it variable coefficients} (see \cite{hy3}).

In this paper we are mainly interested in the practical implementation of the {\it discontinuous} Galerkin neural networks for the {\it nonhomogeneous} governing equations with {\it variable coefficients}. We propose a framework of discontinuous Galerkin neural network method only with a single hidden layer for approximating the quadratic loss functional as in the least square method with $h-$refinement. The desired neural network approximate solution is recursively supplemented by solving a sequence of quasi-minimization problems associated with the underlying loss functionals and the adaptively augmented discontinuous neural network sets without the assumption on the boundedness of the neural network parameters. The convergence results of the proposed method are established. For convenience, we call the method as DGNN method.

For the case of wave equations with constant coefficients, we propose a discontinuous Galerkin Trefftz neural network discretization (DGTNN) with a single hidden layer, which further reduces the computational costs and storage requirements. Moreover, we propose an effective strategy for computing the initial value, which defines a series of nonhomogeneous local problems on auxiliary smooth subdomains by the spectral element method and chooses the associated local spectral elements as the initial value for the DGTNN algorithm. For the case of wave equations with variable coefficients, we propose a discontinuous neural network discretization only with two hidden layers, which employs general activation functions and achieves the given accuracy with relatively smaller computational costs. Particularly, we design a template for initializing nonlinear weights based on the considered equations. Numerous numerical experiments confirm the effectiveness of the DGNN method. Specially, the proposed DGNN method reduces the relative $L^2$ error by at least one order of magnitude under very economical conditions, compared to the existing PINN method and its variations.

The new DGNN method has merits in comparison to most existing methods including the PINN method, its variations and {\it continuous} Galerkin neural networks: (i) the idea upon which it is based is more intuitive and easier to understand even for nonspecialists; (ii) the employed trial and test spaces are the same and consist of discontinuous neural network sets defined {\it locally} over the finite element mesh; (iii) the interface residual losses are further defined by introducing jump operators on the set of all interior faces; and (iv) the template for selecting the initial values of nonlinear weights based on the considered equations determines that the DGNN algorithm is driven by the model.

The paper is organized as follows: in section 2, we give mathematical setup, linear models and discontinuous neural networks.
In section 3, we describe the minimization of a new relaxed loss functional and the associated variational formulation. In section 4, we propose a discontinuous Galerkin neural network iterative method. In section 5, we propose a discontinuous Galerkin Trefftz neural network discretization with a single hidden layer for the case of constant coefficients. In section 6, for the case of variable coefficients, we propose a discontinuous neural network discretization only with two hidden layers. Finally, in section 7, we report some numerical results to confirm the effectiveness of the proposed method.

\section{Mathematical setup, linear models and discontinuous neural networks}
We consider the initial boundary value problem posed on
 a space-time domain $Q=\Omega\times I$, where $\Omega \subset \mathbb{R}^d~ (d\in \mathbb{N})$ is an open bounded Lipschitz polytope and $I=(0,T),T>0$.
The boundary of $\Omega$ denoted by $\Gamma$. We emphasize that the proposed method can certainly be applied to time-harmonic problems.

Let the time domain $(0,T)$ be divided into $N_T\in \mathbb{N}$ intervals $I_p (1\leq p \leq N_T)$ composing a partition ${\cal T}_{h_t}^t$, with
$$I_p=(t_{p-1},t_p),~~h_p=t_p-t_{p-1}=|I_p|, ~~h_t=\mathop\text{max}\limits_{1\leq p \leq N_T}h_p.$$
For each $1\leq p \leq N_T$, we introduce a same polygonal finite element mesh ${\cal T}_{h_{{\bf x}}}^{\bf x}=\{\Omega_q\}_{q=1}^{N_S}$ of the spatial domain $\Omega$ with
$$h_{ \Omega_q }=\text{diam} ~\Omega_q, ~~h_{{\bf x}}=\mathop\text{max}\limits_{\Omega_q\in {\cal T}_{h_{\bf x}}^x}h_{\Omega_q}.$$
Then the space-time domain $Q=\Omega\times(0,T)$ can be partitioned with a finite element mesh ${\cal T}_h$ given by
$${\cal T}_h=\{ Q_k = \Omega_q\times I_p, ~\Omega_q\in {\cal T}_{h_{{\bf x}}}^{\bf x}, ~1\leq p \leq N_T, ~k=1,\cdots,N \}.$$
Here ${\cal T}_h$ is a tensor product mesh.

We denote by $\mathcal{F}_h = \bigcup\limits_{K\in{\cal T}_h}\partial K$ the
skeleton of the mesh.  Set the initial skeleton $\mathcal{F}_h^0 = {\cal T}_{h_{{\bf x}}}^{\bf x} \times \{t=0\}$, the time termination skeleton $\mathcal{F}_h^T={\cal T}_{h_{{\bf x}}}^{\bf x} \times \{t=T\}$,
the boundary skeleton $\mathcal{F}_h^{\text{B}}=\mathcal{F}_h\bigcap (\partial\Omega\times {\cal T}_{h_t}^t)$ and the union of the internal faces
$\mathcal{F}_h^{\text{I}}= \mathcal{F}_h \backslash
(\mathcal{F}_h^{\text{B}}\cup\mathcal{F}_h^0\cup\mathcal{F}_h^T)$.

We denote the outward-pointing unit normal vector on $\partial K$ by $({\bf n}_K^{\bf x},{ n}_K^t)$ for each $K\in{\cal T}_h$.
On an internal face $F=\partial K_1\bigcap \partial K_2$, by convention, we choose ${\bf n}_F^t>0$, which means that the unit normal vector $({\bf n}_F^{\bf x},{\bf n}_F^t)$ points towards future time. ${\bf n}^x_{\Omega}$ is an outward-pointing unit normal vector on $\partial\Omega$.

Let $(V({\cal T}_h),||\cdot||_V)$, $(Y,||\cdot||_Y)$, $(Z,||\cdot||_Z)$ and $(X,||\cdot||_X)$ be Banach spaces.
 Suppose that $\mathcal{A}: V({\cal T}_h)\rightarrow Y$, $\mathcal{B}: V({\cal T}_h)\rightarrow Z$ and $\mathcal{I}: V({\cal T}_h)\rightarrow X$ are bounded linear operators.

Consider an initial boundary value problem (IBVP):
\begin{eqnarray}
\left\{\begin{array}{ll} \mathcal{A}u = f & \text{in}\quad Q,\\
\mathcal{B}u=g& \text{on}\quad\Gamma\times I, \\
\mathcal{I} u(t=0) = q & \text{on}\quad\Omega,
\end{array}\right.
\label{geneform}
\end{eqnarray}
where $f, g, q$ are the given source data. Throughout the paper, we assume that the problem (\ref{geneform}) is well-posed.

{ \bf assumption 1: } %\label{assum}
There exists a unique solution $u\in V({\cal T}_h)$ to the problem (\ref{geneform}).

Suppose that  $\mathcal{C}_i (i=1,\cdots, r): V({\cal T}_h)\rightarrow W_i$ are bounded linear operators related to jump quantities on the set of all interior faces $\mathcal{F}_h^{\text{I}}$, where $(W_i,||\cdot||_{W_i})$ be Banach space. The expressions for jump quantities $\{\mathcal{C}_i\}$ are determined based on the model case by case.

For $Q_k\in {\cal T}_h$, set $u|_{Q_k}=u_k$. Then the reference problem (\ref{geneform}) to be solved consists in finding the local solution $u\in V({\cal T}_h)$ such that

\begin{eqnarray}
\left\{\begin{array}{ll}
\mathcal{A}u_k = f & \text{in}\quad Q_k ~(k=1,\cdots,N),\\
\mathcal{B}u = g &\text{over}\quad \mathcal{F}_h^{\text{B}},\\
\mathcal{I} u(\cdot,0) = q & \text{over}\quad \mathcal{F}_h^0, %\partial Q_k \cap
\end{array}\right. %\quad\quad(k=1,\cdots,N)
\label{helm2}
\end{eqnarray}
and
\be \label{interfacecontinu}
\mathcal{C}_i (u) = 0  \quad \text{over} \quad \mathcal{F}_h^{\text{I}}, \quad i=1,\cdots, r.
\en

Finally, we recall some standard DG notation. Let $w$ and ${\bm \tau}$ be a piecewise smooth function and vector field on ${\cal T}_h$, respectively.
On $F=\partial K_1\bigcap\partial K_2$, we define
\begin{eqnarray}
 & \text{space normal jumps:} ~\llbracket {w} \rrbracket_{\bf N} :=
w_{|K_1}{\bf n}^{\bf x}_{K_1} + w_{|K_2}{\bf n}^{\bf x}_{K_2}, ~\llbracket {\bm \tau} \rrbracket_{\bf N} = {\bm \tau}_{|K_1}\cdot {\bf n}^{\bf x}_{K_1} + {\bm \tau}_{|K_2}\cdot {\bf n}^{\bf x}_{K_2},
\cr &
\text{space tangential jumps:} ~\llbracket {\bm \tau} \rrbracket_{\bf T} = {\bm \tau}_{|K_1} \times {\bf n}^{\bf x}_{K_1} + {\bm \tau}_{|K_2} \times {\bf n}^{\bf x}_{K_2},
\cr &
\text{time full jumps:} ~\llbracket w\rrbracket_t:= w_{|K_1}{ n}^t_{K_1} + w_{|K_2}{ n}^t_{K_2} = (w^- - w^+) { n}^t_F,
\cr &
\text{time full jumps:} ~\llbracket {\bm \tau} \rrbracket_t:= {\bm \tau}_{|K_1}{ n}^t_{K_1} + {\bm \tau}_{|K_2}{ n}^t_{K_2}=({\bm \tau}^- - {\bm \tau}^+){ n}^t_F.
\end{eqnarray}
 Here $w^-$ and $w^+$ denote the traces of the function $w$ from the adjacent elements at lower and higher times, respectively, and similarly for ${\bm \tau}^{\pm}$.

\subsection{Discontinuous neural networks}
A general feed-forward discontinuous neural network consists of a single hidden layer of $n\in \mathbb{R}$ neurons on each element $Q_k\in {\cal T}_h$, defining a function $\varphi^{NN}: \mathbb{R}^{d+1} \rightarrow \mathbb{C}$ as follows:
\be \label{nn_app}
u_{NN}({\bm x}; \Theta)|_{Q_k} = \sum_{j=1}^n c_j^{(k)} \sigma( {\bf W}_j^{(k)} \cdot [{\bm x};t] + b_j^{(k)}), \quad \forall Q_k \in {\cal T}_h,
\en
where $n$ is called as the width of the discontinuous network; ${\bf W}_j^{(k)}\in \mathbb{C}^{d+1}, b_j^{(k)}\in \mathbb{C}$ are elementwise nonlinear parameters; $c_j^{(k)} \in \mathbb{C}$ are elementwise linear parameters; and $\sigma: \mathbb{C}\rightarrow \mathbb{C}$ is a bounded, elementwise smooth activation function. In this paper, we choose the sigmoid activation  function $\sigma(s)=\frac{1}{1+e^{-\xi}}, \xi\in \mathbb{C}$ in  general.

 For simplicity, we shall use the elementwise notation ${\bf W}|_{Q_k}={\bf W}^{(k)}=[{\bf W}_1^{(k)}, \cdots, {\bf W}_n^{(k)}]\in \mathbb{C}^{(d+1)\times n}$ referred as weights, ${\bf b}|_{Q_k}={\bf b}^{(k)}= [b_1^{(k)}, \cdots, b_n^{(k)}]\in \mathbb{C}^{n}$ referred as biases, and ${\bf c}|_{Q_k}={\bf c}^{(k)}=[c_1^{(k)}, \cdots, c_n^{(k)}]^T\in \mathbb{C}^n$. A set of nonlinear and linear parameters over a collection of all elements is denoted by $\Theta=\{{\bf W},{\bf b},{\bf c}\}$. In particular, a set of nonlinear parameters defined on every element is collectively denoted by $\Phi=\{{\bf W},{\bf b}\}$.

Denote $V_n({\mathcal T}_h)$ as the set of all functions of the form (\ref{nn_app}), which is defined by
\be \label{genevsign}
V_n({\mathcal T}_h) := \bigg\{v: ~v|_{Q_k}=\sum_{j=1}^n c_j^{(k)} \sigma( {\bf W}_j^{(k)} \cdot [{\bm x};t] + b_j^{(k)}), %~c_j \in \mathbb{C}, ~W_j \in \mathbb{R}^d, ~x\in\Omega_k,
~Q_k\in {\cal T}_h  \bigg\}.
\en

By the classical approximate procedure \cite{HO}, we can obtain the following universal approximation property of the neural network.

\begin{theorem} \label{geneaninappro}
Suppose that $s\geq 2$, and that $Q\subset \mathbb{R}^{d+1}$ is compact, then $V_n({\mathcal T}_h)$ is dense in the Sobolev space
\be\nonumber
H^{s}({\cal T}_h) := \bigg\{v\in L^2(Q): D^{\alpha}v\in  L^2(Q_k) ~~\forall |\alpha|\leq s, ~Q_k\in {\cal T}_h   \bigg\}.
\en
\end{theorem}

In particular, for any given function
$f \in  H^{s}({\cal T}_h)$ and $\tau  > 0$, there exists $n(\tau , f) \in  \mathbb{N}$  and $\tilde{f} \in  V_{n(\tau ,f )}({\cal T}_{h})$ such that $| | f-\tilde{f}| |_{H^{s}({\cal T}_h)} < \tau$.

However, neural networks exhibit several undesirable properties as well. In particular,
$V_n({\mathcal T}_h)$ is not a vector space since it is not closed under addition when the elementwise nonlinear parameters are not given. But, for simplicity of exposition, we still call $V_n({\mathcal T}_h)$ as a discontinuous neural network space.

The above single hidden layer neural network structures can be extended to multi hidden layer scenarios, especially the case of two hidden layers; see section \ref{twolayer_sec} below.

\section{Minimization of a new relaxed loss functional and associated variational formulation}
A field $u\in V({\cal T}_h)$ is the solution of (\ref{geneform}) if and only if the conditions (\ref{helm2}) and (\ref{interfacecontinu}) are satisfied. Based on this observation, we define a quadratic functional

\beq \label{funct}
& J(v) = \lambda_0\sum_{k=1}^N \int _{Q_k} |\mathcal{A}v_k - f|^2 dV  +
\lambda_1  \int_{\mathcal{F}_h^B} |\mathcal{B}v - g|^2 dS
\cr
& + \lambda_2 \int_{\mathcal{F}_h^0} |\mathcal{I}v(\cdot,0) - q|^2 dS
+ \sum_{i=1}^r \lambda_{i+2}\int_{\mathcal{F}_h^{\text{I}}}
| \mathcal{C}_i (v) |^2 dS, ~~\forall v\in V({\mathcal T}_h),
\eq
where $\{\lambda_i, i=0,\cdots,r+2\}$ are given positive numbers. The introduction of the relaxation parameters $\{\lambda_i\}$ aims to balance the loss terms on the elements, boundaries, initial time and interfaces in the functional $J$. For simplicity, we denote by the PDE residual loss $\mathcal{L}_R(\Theta) = \sum_{k=1}^N \int _{Q_k} |\mathcal{A}u_k - f|^2 dV$, the boundary residual loss $\mathcal{L}_B(\Theta) = \int_{\mathcal{F}_h^B} |\mathcal{B}v - g|^2 dS$, the initial residual loss $\mathcal{L}_0(\Theta) = \int_{\mathcal{F}_h^0} |\mathcal{I}v(\cdot,0) - q|^2 dS$ and the interface residual losses $\mathcal{L}_{I,i}(\Theta) = \int_{\mathcal{F}_h^{\text{I}}}|\mathcal{C}_i (v)|^2 dS, i=1,\cdots, r$.

When the loss $J(v)$ is very small, the field $v$ should be an approximate solution of (\ref{helm2}) and (\ref{interfacecontinu}). According to the above discussion, consider the minimization problem: find $u\in {V}({\cal T}_h)$ such that
\begin{equation} \label{minmum}
J({u})=\min\limits_{{v}\in {V}({\cal T}_h)}~J(v).
\end{equation}

The variational problem associated with the minimization problem (\ref{minmum}) can be expressed as follows: find $u \in V({\mathcal T}_h)$ such that
\begin{eqnarray}
 &
\lambda_0 \sum_{k=1}^N\int _{Q_k}\big(\mathcal{A}u_k - f\big)\cdot\overline{\mathcal{A} v_k}dV
+
 \lambda_1 \int_{\mathcal{F}_h^{\text{B}}} \big(\mathcal{B} u-g\big)\cdot\overline{ \mathcal{B} v}dS \cr
& \quad~
 + \lambda_2 \int_{\mathcal{F}_h^0} \big( \mathcal{I}u(\cdot,0) - q \big) \cdot\overline{  \mathcal{I}v(\cdot,0) } dS
+ \sum_{i=1}^r \lambda_{i+2} \int_{\mathcal{F}_h^{\text{I}}}
\mathcal{C}_i (u)  \cdot \overline{\mathcal{C}_i (v)  } dS=0, ~~\forall v\in V({\mathcal T}_h).
\label{va1}
\end{eqnarray}
Define the sesquilinear form $a(\cdot, \cdot)$ by
\beq \label{sesqu_ele}
&a(u,v) = \lambda_0 \sum_{k=1}^N\int _{Q_k} \mathcal{A}u_k \cdot\overline{\mathcal{A} v_k}dV
 + \lambda_1 \int_{\mathcal{F}_h^{\text{B}}} \mathcal{B} u \cdot\overline{ \mathcal{B} v}dS \cr
 &
 + \lambda_2 \int_{\mathcal{F}_h^0} \mathcal{I}u(\cdot,0) \cdot\overline{  \mathcal{I}v(\cdot,0) } dS
  + \sum_{i=1}^r \lambda_{i+2} \int_{\mathcal{F}_h^{\text{I}}}
 \mathcal{C}_i (u)  \cdot \overline{ \mathcal{C}_i (v) } dS, ~\forall u, v\in V({\mathcal T}_h),
\eq
and semilinear form $L(v)$ by
\be \label{similinear}
L(v)= \lambda_0 \sum_{k=1}^N\int _{Q_k} f \cdot\overline{\mathcal{A} v_k}dV
+ \lambda_1 \int_{\mathcal{F}_h^{\text{B}}} g \cdot\overline{ \mathcal{B} v}dS
+ \lambda_2 \int_{\mathcal{F}_h^0}  q  \cdot\overline{  \mathcal{I}v(\cdot,0) } dS,
 ~\forall v\in V({\mathcal T}_h).
\en
 Then (\ref{va1}) can be written as the following variational formulation:
 \begin{eqnarray}\label{pwlsvar}
\left\{\begin{array}{ll} \text{Find} \ u\in V({\mathcal T}_h) ~~s.t. \\
a(u,v)=L(v) \quad \forall v\in V({\mathcal T}_h).
\end{array}\right.
\end{eqnarray}

\begin{theorem} \label{uniquetheo}
The reference problem (\ref{geneform}) is equivalent to the proposed variational problem (\ref{pwlsvar}).
\end{theorem}

{\it Proof.} It is clear that the solution of the problem (\ref{geneform}) (i.e., the problems (\ref{helm2}))-(\ref{interfacecontinu}) will automatically satisfy the variational problem (\ref{pwlsvar}). Therefore one needs only to verify uniqueness of the solution of the variational problem
(\ref{pwlsvar}). The verification is standard.

 Assume that there are two solutions $u=(u_1,\cdots,u_N)$ and $u'=(u_1',\cdots,u_N')$ of the variational problem. Let $\tilde{u}=(\tilde{u}_1,\cdots,\tilde{u}_N)$ denote the
difference between the two solutions. Plugging the
difference $\tilde{u}$ into the variational form (\ref{pwlsvar}) leads to
\be \label{gene_pr1}
a(\tilde{u},v)=0 \quad \forall v\in V({\mathcal T}_h).
\en
Taking $v=\tilde{u}$, (\ref{gene_pr1}) yields
\be
a(\tilde{u},\tilde{u})=0.
\label{general_eq22}
\en
It implies that
\begin{eqnarray}
\left\{\begin{array}{ll}
\mathcal{A}u_k = 0 & \text{in}\quad Q_k ~(k=1,\cdots,N),\\
\mathcal{B}u = 0 &\text{over}\quad \mathcal{F}_h^{\text{B}},\\
\mathcal{I} u(\cdot,0) = 0 & \text{over}\quad \mathcal{F}_h^0, %\partial Q_k \cap
\end{array}\right. %\quad\quad(k=1,\cdots,N)
\label{gene_pr2}
\end{eqnarray}
and
\be \label{gene_pr3}
\mathcal{C}_i (u) = 0  \quad \text{over} \quad \mathcal{F}_h^{\text{I}}, \quad i=1,\cdots, r.
\en
Combining these with Assumption 1, we know that the function $\tilde{u}$ vanishes on $Q$, which proves the uniqueness of the solution of (\ref{pwlsvar}).
$\hfill\Box$

It is clear that $a(v,v)\geq 0$. Moreover, from the Theorem \ref{uniquetheo}, we can see that $a(v,v)=0$ for
$v\in { V}({\cal T}_h)$ if only if $v=0$. Thus
$a(\cdot,\cdot)$ is a norm on $ { V}({\cal T}_h)$. For ease of
notation, this norm is denoted by $||\cdot||_{ V}^2$. By the Cauchy-Schwarz inequality, we also have that $a(\cdot,\cdot)$ satisfies the boundedness: $|a(u,v)|\leq ||| u |||~||| v |||, ~\forall u,v \in V({\mathcal T}_h)$, and that $L(v)$ is bounded satisfying $|L(v)| \leq ( ||f||_{L^2(Q)} + ||g||_{L^2(\Gamma\times I)} + ||q||_{L^2(\Omega)}) ~||| v |||$.  Furthermore, a standard application of the Risez representation theorem \cite{FO} implies the existence of a unique $u\in V({\mathcal T}_h)$ satisfying (\ref{pwlsvar}).

\begin{remark}
The natural introduction of the PDE residual loss $\mathcal{L}_R(\Theta)$, is attributed to the existence of the source term of the problem the IBVP (\ref{geneform}).
The key difference between the our proposed PDE residual loss and $hp$-VPINNs' is from the employed same trial and test spaces consisting of discontinuous neural network sets defined {\it locally} over the finite element mesh. In order to balance the associated ``large" gradient part caused by the PDE residual loss, adaptive relaxation parameters are introduced to act on other loss parts.
\end{remark}

\section{A discontinuous Galerkin neural network iterative method}\label{helmdis}
To approximately satisfy the minimization problem (\ref{minmum}), we present an iterative approach for quasi-minimizers of the loss functional that combines neural network optimization with discontinuous Galerkin framework.

First, we need to calculate the Fr$\acute{e}$chet derivatives of the defined quadratic functional. For the defined quadratic functional $J: { V}({\cal T}_h)\rightarrow \mathbb{R}$, the first-order Fr$\acute{e}$chet derivative of $J$ at $v\in { V}({\cal T}_h)$ is denoted by $J'(v)$. It is a linear functional (an element of the dual space $(V({\cal T}_h))^{\ast}$), representing the gradient or first variation of $J$ at $v$. Similarly, the second-order Fr$\acute{e}$chet derivative of $J$ at $v\in { V}({\cal T}_h)$ is represented by $J''(v)$, which is a linear operator mapping $V({\cal T}_h)$ to $(V({\cal T}_h))^{\ast}$. Let $\langle \cdot, \cdot \rangle$ denote the duality pairing between $V({\cal T}_h)$ and $(V({\cal T}_h))^{\ast}$, which generalizes the standard $L^2$ inner product to handle derivatives in weak formulations.
It is easy to check that
\be \nonumber
\langle J'(v), w \rangle = 2\text{Re}\big\{ a(v,w)-L(w) \big\},\quad v,w\in { V}({\cal T}_h)
\en
and
\be \label{secondor}
\langle J''(v) w_1,w_2 \rangle = 2\text{Re} \big\{ a(w_1,w_2) \big\}, \quad v,w_1,w_2\in { V}({\cal T}_h).
\en
Here we can see that the second-order Fr$\acute{e}$chet derivation $j''(v)$ is actually independent of the argument $v$.

\subsection{An overview of the DGNN method} \label{dgnn_algo}

Suppose $u_{r-1}\in { V}({\cal T}_h)$ is the current approximation of (\ref{minmum}) at iteration $r-1$. In particular, $u_0$ is an {\it initial guess}, and its suitable selection will improve the accuracy of the proposed algorithm. %We will specifically propose its calculation method in the next section.

Let $\xi_r\in {V}({\cal T}_h)$ be the solution of
the minimization problem:
\begin{equation} \label{iteraminmum}
J(u_{r-1} + {\xi_r})=\min\limits_{{v}\in {V}({\cal T}_h)}~J(u_{r-1} + v).
\end{equation}
This means we are looking for a correction $\xi_r$ in the space
${V}({\cal T}_h)$ that minimizes $J$ when added to the current iterative solution $u_{r-1}$.
At the minimum, the first derivative of $J$ must vanish. This is expressed as:
\be \label{derizero}
\langle J'(u_{r-1} + \xi_r), v \rangle = 0, ~~\forall v \in {V}({\cal T}_h).
\en
Namely, $\xi_r$ satisfies
\be \label{resivari}
a(\xi_r, v) = L(v) - a(u_{r-1}, v),   ~~\forall v \in {V}({\cal T}_h),
\en
 which indicates that $\xi_r$ is the residual of the current iterative solution $u_{r-1}$ and exactly compensates for the error $u - u_{r-1}$:
\be \label{residefi}
\xi_r = u - u_{r-1}.
\en

Since computing the exact minimizer $\xi_r$ of (\ref{iteraminmum}) is equivalent to solving the original minimization problem (\ref{minmum}) explicitly, we instead construct an approximate minimizer $\xi_r^{\theta} \approx \xi_r$ in an appropriate way.

For a positive integer $n_r$ defined later, when aiming to minimize $J(u_{r-1} + \cdot)$ in the subset $V_{n_r}({\mathcal T}_h) \subset {V}({\cal T}_h)$, a significant complication is that the subset $V_{n_r}({\mathcal T}_h)$ is not a linear subspace, nor is it topologically closed in
${V}({\cal T}_h)$. Hence, even though the infimum in $V_{n_r}({\mathcal T}_h)$ of $J(u_{r-1} + \cdot)$ exists (since
$J$ is bounded below by zero), but a minimizer may not exist within $V_{n_r}({\mathcal T}_h)$.
Therefore, one should not aim to completely minimize $J(u_{r-1} + \cdot)$, but instead use a relaxed notion of quasi-minimization as used by \cite{Karniadakis} (for which the existence of an infimum implies the existence of a quasi-minimizer). The idea is that while exact minimization might not be achievable, one can still find an element in $V_{n_r}({\mathcal T}_h)$ that comes ``sufficiently close" to the infimum. Specifically:

 {\bf Definition 1: }
 Let $\delta_r$ be a sufficiently small positive parameter, and let $V_{n_r}({\mathcal T}_h)$ be a subset of ${V}({\cal T}_h)$. A function $\xi_r^{\theta} \in V_{n_r}({\mathcal T}_h)$ is said to be a quasi-minimizer of
 $J(u_{r-1} + \cdot)$ associated with $\delta_r$ if the following inequality holds:
 \be \label{iteradisminmum}
 J(u_{r-1} + \xi_r^{\theta}) \leq \mathop{\text{inf}}\limits_{v_r^{\theta}\in V_{n_r}({\mathcal T}_h)} J(u_{r-1} + v_r^{\theta})+\delta_r.
 \en
 %\end{definition}
It is clear that the value of the functional $J$ evaluated at the quasi-minimizer $\xi_r^{\theta}$ provides a good approximation to the true infimum of $J$.

After computing $\xi_r^{\theta}$, we use this correction term to update the previous approximation $u_{r-1}$, and further get a better approximation $u_r$ to $u$ by setting
\be \label{updateappro}
u_r = u_{r-1} + \xi_r^{\theta}.
\en

We found that it is not feasible to directly employ optimization algorithms to compute discontinuous neural network parameters $\Theta$ for the quasi-minimization problem (\ref{iteradisminmum}), which will lead to unexpected divergence of the algorithms. Instead, the quasi-minimization problem (\ref{iteradisminmum}) is decomposed into two alternating subproblems: optimization of nonlinear parameters and calculation of activation coefficients involving Galerkin least squares type finite element methods. By fixing a set of parameters and optimizing the other (and vice versa), the algorithm iteratively refines the approximation. Besides, adaptive strategy for updating relaxation balancing factors is adopted (see \cite{wangsifan}). {\bf Algorithm 4.1} summarizes this recursive approach.

\vspace{0.2cm}
{\bf Algorithm 4.1} (Alternating training of the quasi-minimization problem (\ref{iteradisminmum}))

Consider the quasi-minimization problem (\ref{iteradisminmum}) with known $u_{r-1}$ and $n_r$. Give initial nonlinear parameters $\Phi^0$ and a tolerance $\rho>0$ for convergence in the iterative process. For known nonlinear parameters $\Phi^l$ with $l\geq 0$,
let $V_{\Phi^l}({\mathcal T}_h)$ denote the space $V_{n_r}({\mathcal T}_h)$ with the given nonlinear parameters $\Phi^l$.

{\bf Step 1} Solve the following minimization problem by the standard Galerkin least squares method
\be \label{generaldglsq}
 J(u_{r-1} + \xi^{l}_r)=\min\limits_{{v}\in V_{\Phi^{l}}({\mathcal T}_h)}~J(u_{r-1} + v)
 \en
to get the activation coefficients ${\bf c}^l$, which are the coordinate vector of the correction term ${\xi^l_r}$ under the activation basis functions with the nonlinear parameters $\Phi^l$;

{\bf Step 2} Update the relaxation parameters by $\lambda_i=(1-\beta)\lambda_i+\beta\hat{\lambda_i}, i\geq 1$, where the recommended hyperparameter $\beta$ is set to be $0.1$, and
\begin{equation} \label{relaxupdate}
\hat{\lambda}_1=\lambda_0\frac{\max_{\Theta_l}{\{|\bigtriangledown_\Theta \mathcal{L}_R(\Theta_l)|\}}}{\overline{\overline{|\bigtriangledown_\Theta(\mathcal{L}_B(\Theta_l))|}}}, ~\hat{\lambda}_2=\lambda_0\frac{\max_{\Theta_l}{\{|\bigtriangledown_\Theta \mathcal{L}_R(\Theta_l)|\}}}{\overline{\overline{|\bigtriangledown_\Theta(\mathcal{L}_0(\Theta_l))|}}}, ~
\hat{\lambda}_{i+2}=\lambda_0\frac{\max_{\Theta_l}{\{|\bigtriangledown_\Theta \mathcal{L}_{R}(\Theta_l)|\}}}{\overline{\overline{|\bigtriangledown_\Theta(\mathcal{L}_{I,i}(\Theta_l))|}}}.
\end{equation}
Here, $\Theta_l$ denotes the values of the parameters $\Theta$ at the $l-$th iteration. $|\cdot|$ denotes the element-wise absolute value, and $\overline{\overline{|\cdot(\Theta_l)|}}$ denotes the mean of $|\cdot|$ at $\Theta_l$.

{\bf Step 3} Solve the following quasi-minimization problem using a gradient descent algorithm,
\be \label{generaldggs}
 J(u_{r-1} + \eta^{l+1}_r) \approx \inf\limits_{{v}\in V_{c^l}({\mathcal T}_h)}~J(u_{r-1} + v),
 \en
 to update the nonlinear parameters $\Phi^{l+1}$ (weights and biases of the activation basis functions) defining the function $\eta^{l+1}_r$,
  where $V_{{\bf c}^l}({\mathcal T}_h)$ is the set $V_{n_r}({\mathcal T}_h)$ with the given activation coefficients ${\bf c}^l$.

If $||\Phi^{l+1} - \Phi^{l}||_{\infty} < \rho$, then the iteration is terminated; %\textcolor{red}{}

else $l=l+1$ and go to {\bf Step 1}. \\

\begin{remark}
The minimization problem (\ref{generaldglsq}) reduces to a linear problem for ${\bf c}^l$ once the nonlinear parameters $\Phi^l$ are typically fixed. Specially, the Galerkin least squares method minimizes $J(u_{r-1} + v)$ by requiring the sum of the current iterative solution $u_{r-1}$ and its residual $\xi^{l}_r$ to be solved to be orthogonal to the discrete test space $V_{\Phi^{l}}({\mathcal T}_h)$ in the definition of the sesquilinear form $a(\cdot, \cdot)$. This is equivalent to solving
\be \label{dresivari}
a(\xi^{l}_r, v) = L(v) - a(u_{r-1}, v),   ~~\forall v \in V_{\Phi^{l}}({\mathcal T}_h).
\en
\end{remark}

 \begin{remark}
 In order to compute optima of objective functions (\ref{generaldggs}), we need to compute the gradient of the objective function $J(u_{r-1} + \eta_r)$ with respect to the nonlinear parameters $\Phi$:
\be
\nabla_{\Phi} J(u_{r-1} + \eta_r),
\en
and perform gradient descent steps for the nonlinear parameters $\Phi$ with learning rate $\alpha$:
\be \label{gradientIMPLE}
\Phi^{l,k+1} = \Phi^{l,k} - \alpha \nabla_{\Phi} J(u_{r-1} + \eta_r^{l,k}),
\en
where
$k$ is the training iteration index, $\eta_r^{l,k}$ is defined by the current parameters $\Phi^{l,k}$ (Note that $\Phi^{l,0}=\Phi^l$). The process continues until convergence to $\Phi^{l+1}$.
\end{remark}

 A pseudo {\bf Algorithm 1} summarizes this recursive approach.

 \begin{algorithm}[H]
\caption{Discontinuous Galerkin Neural Network Iteration}

\quad Set $u_0$ and tolerance $tol > 0$.

\quad {\bf for} $r=1 : \text{maxit}$ {\bf do}

\quad\quad  Given $u_{r-1}$, $n_r$ and a small parameter $\rho>0$.

\quad\quad  Initialize hidden parameters $\Phi^0\in \mathbb{C}^{n_r}$ case by case.

\quad\quad  {\bf for} $l = 0 : \text{traincount}$

\quad \quad \quad Compute ${\bf c}^l$ by the standard Galerkin least squares
method.

\quad \quad \quad Update relaxation parameters $\{\lambda_i\}$.

\quad \quad  \quad Train nonlinear parameters $\Phi^{l+1}$.

\quad \quad  \quad {\bf if } $||\Phi^{l+1} - \Phi^{l}||_{\infty} < \rho$

\quad \quad  \quad \quad return $\xi_r^l$ and $J(u_{r-1}+\xi_r^l)$

\quad \quad  \quad {\bf end if}

\quad\quad  {\bf end for}

\quad \quad $u_{r} = u_{r-1} +\xi_r^l$

\quad \quad $r = r+1$

\quad \quad {\bf if } $J(u_{r-1}) < tol$ or relative $L^2$-norm error $< tol$

\quad \quad \quad break

\quad \quad {\bf end if}

\quad  {\bf end for}

\quad Return $N=r-1$ and $u_N$.
\end{algorithm}

\begin{figure}[H]
%\vspace{-2cm}thb
\begin{center}
\begin{tabular}{c}
\epsfxsize=0.8\textwidth\epsffile{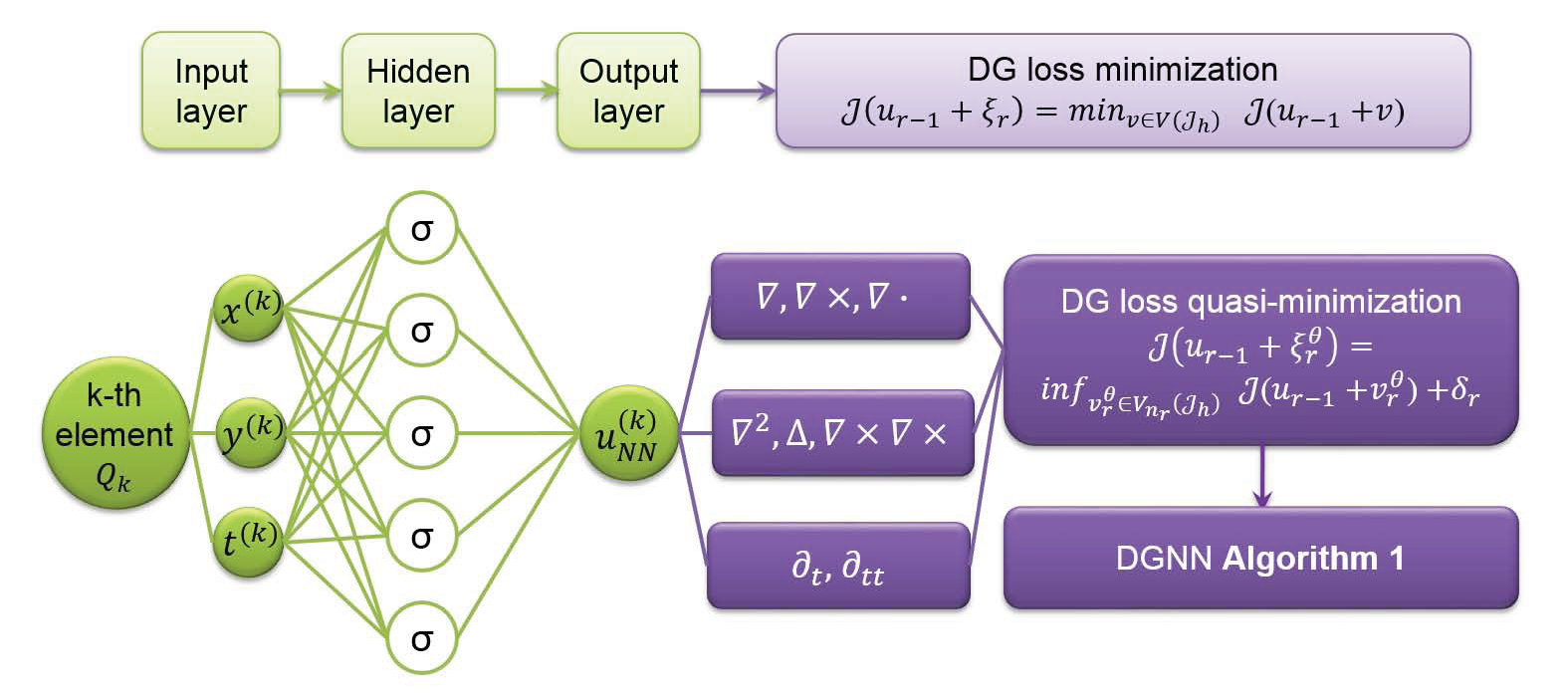}\\
\end{tabular}
\end{center}
 \caption{ A schematic for the discontinuous Galerkin Neural network method with a single hidden layer. }
\label{1hlnn}
\end{figure}

Figure \ref{1hlnn} shows a schematic for our DGNN method with a single hidden layer. % in a DGNN iterative method.

\subsection{The convergence analysis}

The next result addresses the convergence of the proposed recursive algorithm.
We mainly used the approximation property of general activation functions, as shown in Theorem \ref{geneaninappro}, and the special property of the new quadratic functional, that is, the second derivative $J''(v)$ is constant independent of the argument $v\in
V({\mathcal T}_h)$.

\begin{theorem} \label{helmtheom}
Assume that, for $j=1,\cdots,r$, $|||u-u_{j-1}|||>0$,~$0<\tau_j < 1$ and $\delta_j=\frac{\tau^2_j}{2}~ |||u-u_{j-1}|||^2$. There exists $n(\tau_j, u_{j-1})\in \mathbb{N}$ such that, when $n_j\geq n(\tau_j, u_{j-1})$, the quasi-minimizer $\xi_j^{\theta}\in V_{n_j}({\mathcal T}_h)$ defined in Definition 1 has the approximate estimate:
\be \label{redualappro}
||| \xi_j^{\theta} - \xi_j ||| < \tau_j ~|||u-u_{j-1}|||.
\en
Further, the resulting approximate solution $u_r$ has the error estimate:
\be
|||u-u_r||| < |||u-u_0|||~\prod_{j=1}^r\tau_j.
\en
\end{theorem}

{\it Proof.}
Define $j(v)= J(u_{r-1} + v)$. Then by (\ref{derizero}), we have $j'(\xi_r)=0$.
\iffalse
Similarly with (\ref{secondor}), we have
\be \label{dissecond}
 \langle j''(v)w_1,~w_2\rangle= 2\text{Re}\big\{ a(w_1,~w_2) \big\},\quad \forall v,w_1,w_2 \in V({\cal T}_h).
 \en
\fi
By the first-order Taylor formula with Lagrange type remainder, there exists a function $\xi^{\ast}_r\in
V({\mathcal T}_h)$ such that %(see \cite)
\beq \label{tay1}
 & j(\xi_r^{\theta}) - j(\xi_r)  =  \langle j'(\xi_r), \xi_r^{\theta}-\xi_r \rangle  + \frac{1}{2} \langle j''(\xi^{\ast}_r)(\xi_r^{\theta}-\xi_r),\xi_r^{\theta}-\xi_r \rangle \cr
& = \text{Re}\big\{ a(\xi_r^{\theta}-\xi_r,~\xi_r^{\theta}-\xi_r) \big\}
  = a(\xi_r^{\theta}-\xi_r,\xi_r^{\theta}-\xi_r) = ||| \xi_r^{\theta} - \xi_r|||^2.
%& \frac{1}{2}  \langle j''(\xi^{\ast}_r)(\xi_r^{\theta}-\xi_r),\xi_r^{\theta}-\xi_r \rangle = a(\xi_r^{\theta}-\xi_r,\xi_r^{\theta}-\xi_r)  = ||| \xi_r^{\theta} - \xi_r|||^2,
\eq
Here the formula (\ref{secondor}) has been used in the second equality.
Similarly to (\ref{tay1}), we get
\be\label{tay2}
j(v_r^{\theta}) - j(\xi_r) =  |||v_r^{\theta}-\xi_r |||^2, \quad~\forall v_r^{\theta}\in V_{n_r}({\mathcal T}_h).
\en
Using quasi-minimization property of $j(\xi_r^{\theta})$, we have
\be \label{quaimini}
 j(\xi_r^{\theta}) - j(\xi_r) \leq \mathop{\text{inf}}\limits_{v_r^{\theta}\in V_{n_r}({\mathcal T}_h)} j(v_r^{\theta}) + \delta_r - j(\xi_r).
 \en
 %\leq j(v_r^{\theta}) - j(\xi_r) + \delta_r.$$
Combining (\ref{tay1}), (\ref{tay2}) with (\ref{quaimini}), we get the quasi-optimal estimate
\be \label{bestappro}
|||\xi_r - \xi_r^{\theta} |||^2 \leq \mathop{\text{inf}}\limits_{v_r^{\theta}\in V_{n_r}({\mathcal T}_h)}|||\xi_r - v_r^{\theta} |||^2 + \delta_r.
\en
Moreover, by the universal approximation property Theorem \ref{geneaninappro}, there exists $n(\tau_r, u_{r-1})\in \mathbb{N}$ and $\tilde{\xi}_r^{\theta}\in  V_{n(\tau_r, u_{r-1})}({\mathcal T}_h)$ such that
\be \label{interappro}
|||\xi_r-\tilde{\xi_r^{\theta}}||| < \frac{\tau_r}{\sqrt{2}}~ |||u-u_{r-1}|||.
\en
Assume $n_r\geq n(\tau_r, u_{r-1})$ and choose $v_r^{\theta}=\tilde{\xi_r^{\theta}}$ in (\ref{bestappro}). Pluging (\ref{interappro}) into (\ref{bestappro}) and
noting that $\delta_r= \frac{\tau^2_r}{2}~ |||u-u_{r-1}|||^2$, we get the desired result (\ref{redualappro}).

By utilizing the iterative scheme (\ref{updateappro}) and the definition (\ref{residefi}) of the residual $\xi_r$, we obtain
\be \label{initierror}
 |||u-u_r||| = ||| u - ( u_{r-1} + \xi_r^{\theta}) ||| = ||| \xi_r -  \xi_r^{\theta} |||.
 \en
Recursively using the approximation estimate (\ref{redualappro}), yields
%can directly the result.
\beq \label{iteresti}
|||u-u_r||| & < & \tau_r |||u-u_{r-1}||| \cr
&& \vdots
 \cr
 & < &  |||u-u_0|||~\prod_{j=1}^r \tau_j.
\eq
$\hfill\Box$

\begin{remark}
Although the source term in nonhomogeneous equations leads to the occurrence of elemental integrals in PDE residual loss functionals, the total loss functional $J$ is still quadratic convex. This, along with the approximation of activation functions guarantees the above textual error estimates. We would like to point out that, the considered nonhomogeneous model determines the convergence speed of the algorithm, which also depends on initial hidden parameters, the initial value computing, approximation capability of activation functions such as sigmoid function, and learning efficiency of the alternating training algorithm proposed. In particular, for the DPWNN method in \cite{yuanhu}, the plane wave activation function used does not have the ability to approximate nonhomogeneous functions, making the DPWNN method ineffective for solving the current model.
\end{remark}

\begin{remark} %\textcolor{red}{initializing hidden parameters}
To enhance the robustness of training of the DGNN method and the accuracy of its predictions,
initial hidden parameters $\Phi$ should be set according to the model case by case, e.g. see section \ref{variahelm} for setting initial values $\Phi$ for two hidden layer neural networks.
\end{remark}

\begin{remark} %\justifying
Adaptive relaxation parameters are introduced by utilizing the back propagated gradient statistics during model trainings,
such that the interplay among different loss terms in composite loss functions is appropriately balanced.
\end{remark}

\begin{remark} %\textcolor{red}{singularity}: ~~~
The proposed DGNN algorithm can introduce the auxiliary variable to reduce the original high-order differential equation to a low-order differential system, and thus efficiently deal with rough solutions such as singularities and sharp changes on interfaces. Moreover, the associated loss functional relaxes the smoothness requirements on the analytic solutions and the discrete spaces induced by activation functions. For further progress, please refer to our future articles.
\end{remark}

\subsection{Discussions on the default initial value $u_0=0$}
To illustrate the impact of the initial value $u_0$ on the algorithm accuracy, we use the following examples to demonstrate the performance of the algorithm when the initial value is set to be {\it zero by default}.

\subsubsection{Possion equations}
Consider the Possion model
\begin{eqnarray}
\left\{\begin{array}{ll} -\Delta u = f & \text{in}\quad
\Omega,\\
u=g& \text{on}\quad\Gamma=\partial\Omega.
\end{array}\right.
\end{eqnarray}
The exact solution of the above problem is defined in the closed form
$$u_{ex}(x,y)= x ~\text{cos} y$$
with $\Omega=[0,1]\times[0,1]$.
The source term $f= x ~\text{cos} y.$

We also denote by $\mathcal{F}_h = \bigcup\limits_{K\in{\cal T}_h}\partial K$ the skeleton of the spacial mesh ${\cal T}_h$.  Set the boundary skeleton $\mathcal{F}_h^{\text{B}}=\mathcal{F}_h\bigcap \partial\Omega$ and the union of the internal faces
$\mathcal{F}_h^{\text{I}}= \mathcal{F}_h \backslash \mathcal{F}_h^{\text{B}}$.

Obviously, the operators $\mathcal{A}, \mathcal{B}, \mathcal{C}_1, \mathcal{C}_2$ are defined by, respectively,
\beq \label{fpossionoperadefin}
& \mathcal{A}= -\Delta ~\text{in}~ \Omega, \quad \mathcal{B} = I ~\text{on}~ \Gamma,
\cr
& \mathcal{C}_1 (u)  = \llbracket u \rrbracket =u_k-u_j, \quad \mathcal{C}_2 (u) = \llbracket \nabla u \rrbracket_{\bf N}= \partial_{\text{\bf n}_k}u_k + \partial_{\text{\bf n}_j}u_j ~\text{on}~ \mathcal{F}_h^{\text{I}},
%\Gamma_{kj},
\eq
where the outer normal derivative is referred to by $\partial_{\bf n}$.

The relaxation parameters $\{\lambda_i\}$ are initialized as follows.
\be \nonumber
\lambda_0=\lambda_3=1, ~\lambda_1=\lambda_2=200\pi.
\en
%\textcolor{red}{
 Weights ${\bf W}^{(k)}\in \mathbb{C}^{n\times 2}$ are initialized to be uniformly distributed.
\be \label{direcini}
{ W}^{(k)}_j = (cos\theta_j, sin\theta_j), ~~~\theta_j=\frac{2\pi}{n}j, j=1,\cdots,n.
\en
Set biases ${\bf b}^{(k)}={\bf 1}\in \mathbb{C}^{n}$.
%}
Set $h=\frac{1}{16}$ and $tol = 10^{-3}$.

\begin{figure}[H]
%\vspace{-2cm}thb
\begin{center}
\begin{tabular}{cc}
\epsfxsize=0.4\textwidth\epsffile{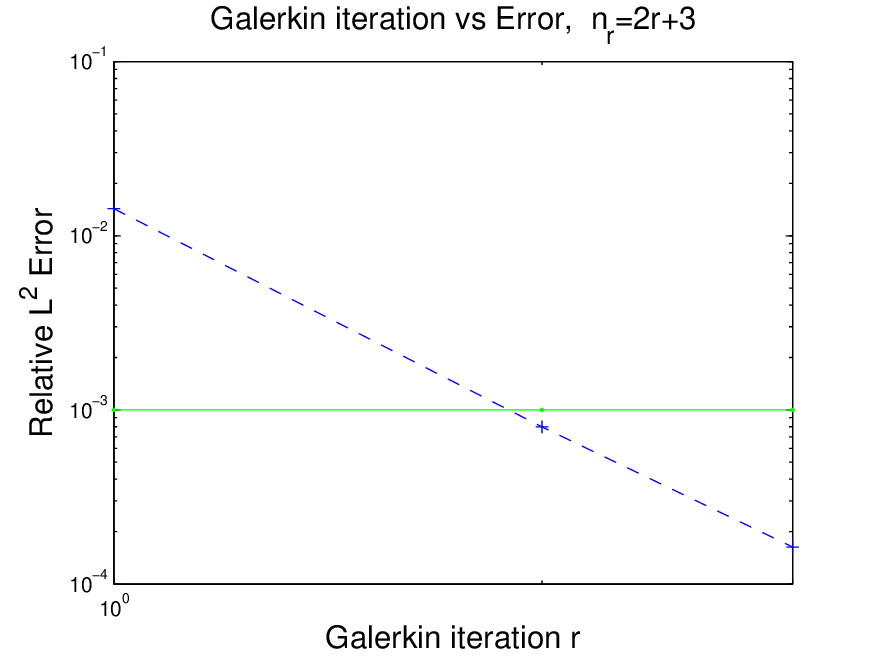}&
\epsfxsize=0.4\textwidth\epsffile{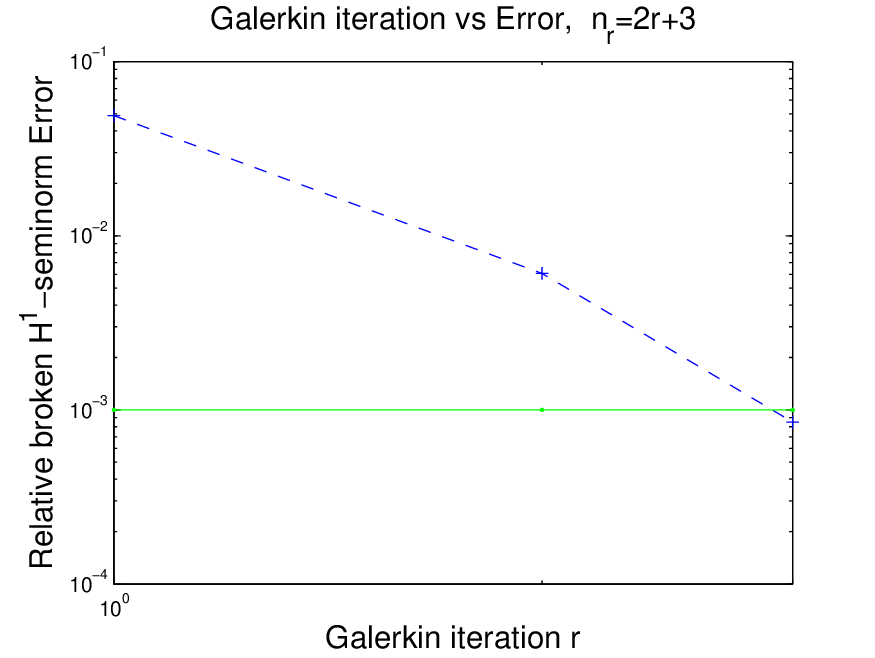}\\
\end{tabular}
\end{center}
 \caption{Displacement of a string (Possion equation in two dimension). (Left) Relative error in the $L^2$-norm at each Galerkin iteration. (Right) Relative error in the broken $H^1$-seminorm at each Galerkin iteration. }
\label{possionnn1}
\end{figure}

Figure \ref{possionnn1} shows the relative errors of $u-u_r$ in the $L^2$-norm and the broken $H^1$-seminorm at each Galerkin iteration, respectively.

Numerical results validate that, when choosing default initial value $u_0=0$, the resulting approximate solutions can reach the given accuracy.

\subsubsection{Helmholtz equations}
Consider the Helmholtz model
\begin{eqnarray}
\label{helmmodel}
\left\{\begin{array}{ll} -\Delta u-\omega^2 \rho u=f& \text{in}\quad
\Omega,\\
(\partial_\text{\bf n} + \text{i}\omega)u=g& \text{on}\quad\Gamma=\partial\Omega.
\end{array}\right.
\end{eqnarray}
%and the angular frequency by $\omega$.}
Set $\rho=1$. The exact solution of the above problem is defined in the closed form
$$u_{ex}(x,y)=\omega x ~\text{cos}y+ y~\text{sin}(\omega x)$$
with $\Omega=[0,1]\times[0,1]$, and $g=({\partial \over
\partial {\bf n}}+i\omega)u_{ex}$.
The source term $f=(1-\omega^2)\omega x ~\text{cos} y.$

Clearly, the operators $\mathcal{A}, \mathcal{B}, \mathcal{C}_1, \mathcal{C}_2$ are defined by, respectively,
\beq \label{foperadefin}
& \mathcal{A}= -\Delta -\omega^2 \rho ~\text{in}~ \Omega,  \quad \mathcal{B} = \partial_\text{\bf n} + \text{i}\omega ~\text{on}~ \Gamma,
\cr
& \mathcal{C}_1 (u)  = \llbracket u \rrbracket , \quad
 \mathcal{C}_2 (u) = \llbracket \nabla u \rrbracket_{\bf N} ~\text{on}~ \mathcal{F}_h^{\text{I}}.
%\Gamma_{kj},
\eq

%Initialize nonlinear parameters as follows.
The relaxation parameters $\{\lambda_i\}$ are initialized as follows.
\be\nonumber
 ~\lambda_0=1, ~\lambda_1=\lambda_3=\omega^2, \lambda_2=\omega^4.
\en
 Weights ${\bf W}^{(k)}\in \mathbb{C}^{n\times 2}$ are initialized to be uniformly distributed.
\be\nonumber
{ W}^{(k)}_j = \text{i}\omega(cos\theta_j, sin\theta_j), ~~~\theta_j=\frac{2\pi}{n}j, j=1,\cdots,n.
\en
Set biases ${\bf b}^{(k)}={\bf 1}\in \mathbb{C}^{n}$. The other parameters are set as follows.
\be \nonumber
%\omega& = & 16\pi, 32\pi, ~ h=\frac{1}{8}, ~n_r = 4r+1, ~tol = 10^{-6}; \cr
 \omega=4\pi, ~h=\frac{1}{20}, n_r = 2r+9, ~tol = 10^{-3}.
\en

\begin{figure}[H]
%\vspace{-2cm}thb
\begin{center}
\begin{tabular}{cc}
\epsfxsize=0.4\textwidth\epsffile{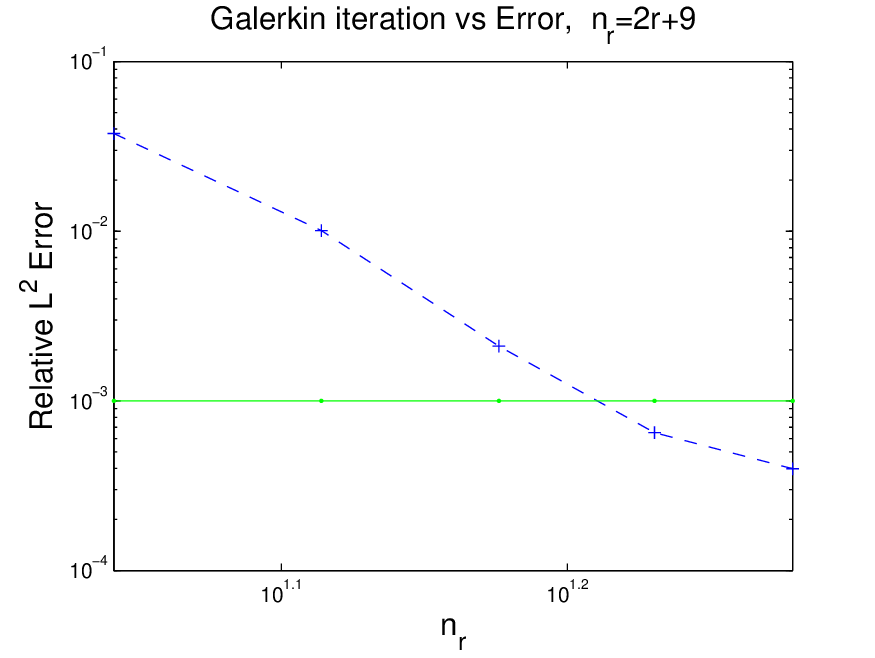}&
\epsfxsize=0.4\textwidth\epsffile{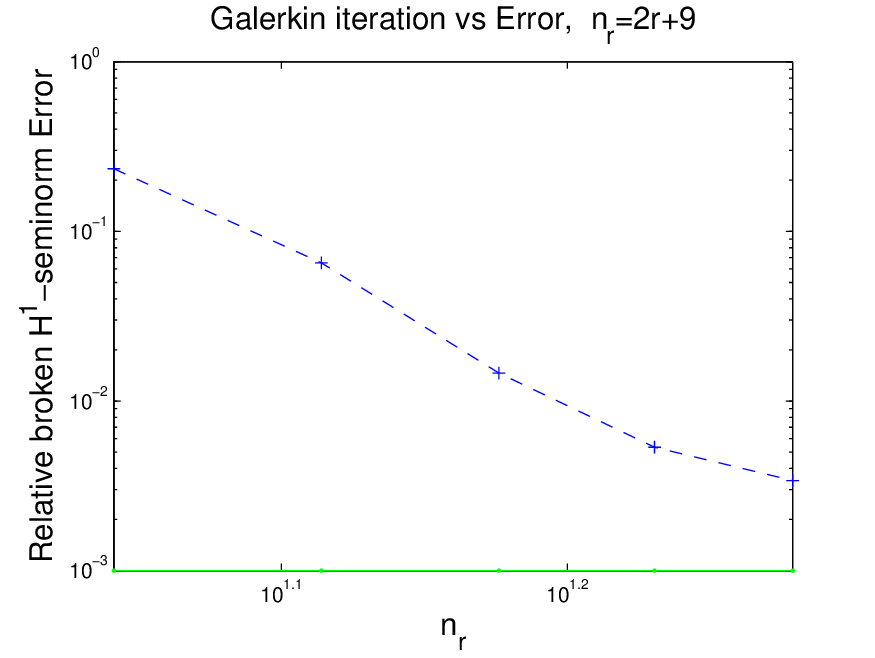}\\
\end{tabular}
\end{center}
 \caption{ Displacement of a string (Helmholtz equation in two dimension). (Left) Relative error in the $L^2$-norm at each Galerkin iteration. (Right) Relative error in the broken $H^1$-seminorm at each Galerkin iteration. }
\label{helmnn1}
\end{figure}

Figure \ref{helmnn1} shows the relative errors of $u-u_r$ in the $L^2$-norm and the broken $H^1$-seminorm at each Galerkin iteration, respectively.

%\textcolor{red}{need to modify}

Numerical results validate that, when choosing the default initial value $u_0=0$ and the single layer neural network architecture with general activation functions,
 to achieve the specified accuracy, the DGNN algorithm requires small step size $h$ and a large number of neurons $n_r$ on each element. It motivates us to propose appropriate strategies for computing initial values and employ neural network structures with some particular activation functions adapted to specific models.

\section{A high-accuracy DGTNN method with a single hidden layer for PDEs with constant coefficients } \label{constmodel}
In this section, we focus on employing a particular neural network structure with the Trefftz activation functions for a dramatic increase in efficiency and accuracy when solving specific types of differential equations with constant coefficients. For this, we need to calculate an appropriate initial value first.

\subsection{An effective strategy for computing the initial value $u_0$} \label{local_spec}
In order to effectively handle the nonhomogeneous model with the source term $f$, as in \cite{hy3}, we
first solve a series of nonhomogeneous local problems on auxiliary smooth subdomains by the spectral element method,
 and choose the associated local spectral elements as the initial value $u_0$ for the original DGNN iterative algorithm.
 Then the DGNN method is applied to the resulting minimization of loss functional on the discrete neural network spaces.

For it, we need to assume that the function $f$ is well defined in a slightly large domain containing $\Omega$ as its subdomain, and then decompose the solution $u$ of the problem \((\ref{geneform})\) into \( u = u^{(1)} + u^{\text{nn}} \), where \(u^{(1)}\) is a particular solution of (\ref{geneform}) with the homogeneous artificial boundary condition, and $u^{\text{nn}}$ needs to make $u$ to be the solution to the minimization problem (\ref{minmum}).% associated with the globally homogeneous neural network sets.

Specifically, for each element $\Omega_q\in \Omega$, let $\Omega^{\ast}_q$ be a fictitious domain that has almost the same size with $\Omega_q$ and contains $\Omega_q$ as its subdomain. Let $u^{(1)}\in L^2(Q)$ be defined as $u^{(1)}\mid_{Q_k}=u^{(1)}_k\mid_{Q_k}$ for each $Q_k$, where $u^{(1)}_k$ satisfies the nonhomogeneous {\it local} governing equation on the fictitious domain $ Q^{\ast}_k = \Omega^{\ast}_q \times I_p$:

\begin{eqnarray}
\left\{\begin{array}{ll}
\mathcal{A}u^{(1)}_k = f & \text{in}\quad Q^{\ast}_k,\\
\mathcal{B}u^{(1)}_k = 0 &\text{on}\quad \partial\Omega^{\ast}_q \times I_p,\\
\mathcal{I} u^{(1)}_k(\cdot, t_{p-1}) = 0 & \text{on}\quad \Omega^{\ast}_q\times \{t=t_{p-1}\},
\end{array}\right.\quad\quad(k=1,\cdots,N).
\label{cha4.n1}
\end{eqnarray}
Clearly, we can find the weak solution \(u^{(1)}_k\) satisfying the corresponding variational problem of (\ref{cha4.n1}), which will be further identified in sections \ref{conshelmdgtnn} and \ref{numerical} by case-by-case.

Let $m$ be a positive integer and $D$ be a bounded and connected domain in $\mathbb{R}^{d+1}$. Denote by $S_m(D)$ the set of polynomials defined on $D$, whose orders are less or equal to $m$. Let \(u_{k,h}^{(1)}\) introduced in section \ref{numerical} denote the approximation of the resulting discrete variational problem of (\ref{cha4.n1}) involved on the local discrete space $S_m(Q^{\ast}_k)$. Define $u^{(1)}_h \in \prod_{k=1}^N S_m(Q_k)$ by $u^{(1)}_h|_{Q_k}=u^{(1)}_{k,h}|_{Q_k}$, and choose $u^{(1)}_h$ as the initial solution $u_0$ for the DGNN iterative algorithm.

It is immediate to see that \( u^{\text{nn}} = u-u^{(1)}\) is uniquely determined by the following minimization problem: find $ u^{\text{nn}} \in {V}({\cal T}_h)$ such that
\begin{equation} \label{minmumimprove}
J({u^{(1)} + u^{\text{nn}}})=\min\limits_{v^{\text{nn}}\in {V}({\cal T}_h)}~J(u^{(1)} + v^{\text{nn}}).
\end{equation}
We seek a discrete quasi-minimizer $\xi_r^{\theta}$ of the residual
 $J( u_{r-1} + \cdot)$ (Note that $u_0=u^{(1)}_h$) associated with $\delta_r$ satisfying:
 \be \label{iteradisminmum_improve}
 J( u_{r-1} + \xi_r^{\theta}) \leq \mathop{\text{inf}}\limits_{v_r^{\theta}\in V_{n_r}({\mathcal T}_h)} J( u_{r-1} + v_r^{\theta})+\delta_r,
 \en
where $V_{n_r}({\mathcal T}_h)$ is the associated discrete neural network set. % with $n_r$ degrees of freedom on each element.

After solving $\xi_r^{\theta}$, we use it to further get a better approximation $u_r$ to $u$ by setting
\be \label{updateappro_improve}
u_r = u_{r-1} + \xi_r^{\theta}.
\en

As in the proof of Theorem \ref{helmtheom}, we can derive the convergence rates below.

\begin{theorem} \label{helmtheom_improve}
Assume that, for $j=1,\cdots,r$, $|||u - u_{j-1}|||>0$,~$0<\tau_j < 1$ and $\delta_j=\frac{\tau^2_j}{2}~ |||u- u_{j-1}|||^2$. Let $n(\tau_j, u_{j-1})\in \mathbb{N}$ and
$\xi_j^{\theta}\in V_{n_j}({\mathcal T}_h)$ with $n_j\geq n(\tau_j, u_{j-1})$ be determined as in Theorem \ref{helmtheom} ($j=1,\cdots,r$).
Then the resulting approximate solution $u_r$ has the error estimate:
\be
\label{impro_conver}
|||u- u_r||| < |||u- u_h^{(1)}|||~\prod_{j=1}^r\tau_j.
\en
\end{theorem}

\begin{remark}
We would like to emphasize that the problems (\ref{cha4.n1}) are
local and independent each other for $k = 1,2,\ldots,N$, so they can be solved in parallel
and the cost is very small. %\textcolor{red}{good initial approximation}
\end{remark}

\subsection{A discontinuous Galerkin Trefftz neural network method} \label{DGTNN}
As we know, one of the main difficulties in the numerical simulation of the propagating problem arises from the usually oscillatory nature of the solutions. Trefftz methods have been designed for finite element discretization of time-harmonic and time-dependent wave equations (for example, the Helmholtz and Maxwell equations) in the last years, two typical Trefftz methods are the plane wave discontinuous Galerkin (PWDG) method \cite{ref21,pwdg,YHESAIM} and the plane wave least squares (PWLS) method \cite{hy2,hy3,peng}. The Trefftz  methods have an important advantage over Lagrange finite elements for discretization of the wave equations: to achieve the same accuracy, relatively smaller number of degrees of freedom are enough in the Trefftz-type methods owing to the particular choice of the basis functions that satisfy the considered PDEs without boundary conditions.
%However the standard Trefftz methods are heavily unstable \cite{ref11,ref21,HMPsur,HGA}, which limits real accuracies of the resulting approximate solutions.
At present, we try to combine the proposed discontinuous Galerkin neural network method with the Trefttz method for the discretization of wave equations with constant coefficients.
To this end, let $\text{T}(Q_k)$ denote the space of the functions which verify the homogeneous governing equation on each element $Q_k$:
\be
 \text{T}(Q_k) = \bigg\{v_k \in L^2(Q_k): ~~\mathcal{A} v_k = 0 \bigg\}.
 \en
Define the following Trefftz space:
\be
\text{T}({\cal T}_h) = \{ v\in L^2(Q): ~v|_{Q_k}\in \text{T}(Q_k) ~~\forall Q_k\in {\cal T}_h  \},
\en
which is also referred to as the globally homogeneous neural network set.
Then for $v^{\text{nn}}\in \text{T}({\cal T}_h)$, the loss functional $J(u^{(1)} + v^{\text{nn}})$ degenerates into a simple form containing only boundary integrals:
 \beq \label{funct2}
 & J(u^{(1)} + v^{\text{nn}})  = \lambda_1  \int_{\mathcal{F}_h^B} |\mathcal{B} (u^{(1)} + v^{\text{nn}}) - g|^2 dS
\cr
& + \lambda_2 \int_{\mathcal{F}_h^0} |\mathcal{I} (u^{(1)} + v^{\text{nn}})(\cdot,0) - q|^2 dS
+ \sum_{i=1}^r \lambda_{i+2}\int_{\mathcal{F}_h^{\text{I}}}
| \mathcal{C}_i (u^{(1)} + v^{\text{nn}}) |^2 dS.
\eq

Consider the minimization problem: find $u^{\text{nn}}\in \text{T}({\cal T}_h)$ such that
\begin{equation} \label{treminmum}
J(u^{(1)} + u^{\text{nn}})=\min\limits_{v^{\text{nn}}\in \text{T}({\cal T}_h)}~J(u^{(1)} + v^{\text{nn}}).
\end{equation}
The variational problem associated with the minimization problem (\ref{treminmum}) can be expressed as follows: find $u^{\text{nn}}\in \text{T}({\cal T}_h)$ such that
 \begin{eqnarray}\label{dgtnnvar}
\left\{\begin{array}{ll} \text{Find} \ u^{\text{nn}} \in \text{T}({\cal T}_h) ~~s.t. \\
a(u^{\text{nn}},v) = L(v) - a(u^{(1)},v) \quad \forall v\in \text{T}({\mathcal T}_h).
\end{array}\right.
\end{eqnarray}

We would like to emphasize that all the elemental integrals included in the variational problem (\ref{va1}), sesquilinear form $a(\cdot, \cdot)$ defined by (\ref{sesqu_ele}) and semilinear form $L(v)$ defined by (\ref{similinear}) vanish in the present situation. Of course, the associated elemental integrals also do not contribute to the gradient descent method in Step 2 of Algorithm 4.2.

A discrete quasi-minimizer $\xi_r^{\theta}$ of the above residuals $J( u_{r-1} + \cdot)$ (Note that $u_0=u^{(1)}_h$) satisfies
 \be \label{iteradisminmum_improvetre}
 J( u_{r-1} + \xi_r^{\theta}) \leq \mathop{\text{inf}}\limits_{v_r^{\theta}\in \text{T}_{n_r}({\mathcal T}_h)} J( u_{r-1} + v_r^{\theta})+\delta_r.
 \en
%(\ref{iteradisminmum_improve}) and can be found in the globally homogeneous neural network sets $\text{T}_n({\mathcal T}_h)$.
Some related discontinuous Trefftz neural network sets $\text{T}_n({\mathcal T}_h)$ with a single hidden layer can be found in sections \ref{conshelmdgtnn} and \ref{numerical}.

\subsection{Helmholtz equations with constant coefficient $\rho=1$}\label{conshelmdgtnn}
Reconsider the Helmholtz model (\ref{helmmodel}). The exact solution of the above problem is defined in the closed form
$$u_{ex}(x,y)=\omega x ~\text{cos}y+ y~\text{sin}(\omega x)$$
with $\Omega=[0,1]\times[0,1]$, and $g=({\partial \over
\partial {\bf n}} + \text{i}\omega)u_{ex}$.
The source term $f=(1-\omega^2)\omega x ~\text{cos} y.$

As introduced in the previous section \ref{local_spec} for computing the initial value $u_0$, the variational problem of the nonhomogeneous local governing equation \((\ref{cha4.n1})\) is to find \(
u^{(1)}_k\in H^1(\Omega^{\ast}_k)\) such that
\begin{equation} \label{eq4.n2} \left\{
\begin{aligned}
     &\int_{\Omega^{\ast}_k}(\nabla u^{(1)}_k\cdot \overline{\nabla{v}_k}  - \omega^2 \rho u^{(1)}_k \overline{v_k} )dV+
\int_{\partial\Omega^{\ast}_k} \text{i}\omega u^{(1)}_k \overline{v_k} dS =
\int_{\Omega^{\ast}_k}f ~\overline{v_k} dV,\\
    & \quad\quad\quad\quad\quad\forall v_k\in
H^1(\Omega^{\ast}_k)~(k = 1,2,\ldots,N).
                          \end{aligned} \right.
                          \end{equation}

We choose local spectral elements of order $m=3$ as the initial value $u_0$ for the original DGTNN iterative algorithm. The discrete variational problems of (\ref{eq4.n2}) are: to find $u^{(1)}_{k,h}\in S_3(\Omega^{\ast}_k)$ such that
\begin{equation} \label{3.new1} \left\{
\begin{aligned}
     &\int_{\Omega^{\ast}_k}(\nabla u^{(1)}_{k,h}\cdot\nabla \overline{v_k}  - \omega^2 \rho u^{(1)}_{k,h} \overline{v_k} )dV +
\int_{\partial\Omega^{\ast}_k} \text{i}\omega u^{(1)}_{k,h} \overline{v_k} dS=
\int_{\Omega^{\ast}_k}f ~\overline{v_k} dV,\\
    & \quad\quad\quad\quad\quad\forall v_k\in
S_3(\Omega^{\ast}_k)~(k = 1,2,\ldots,N).
                          \end{aligned} \right.
                          \end{equation}

 We will choose two types of activation functions to generate neural network spaces. One is the standard sigmoid activation function, which associates with the DGNN method introduced in section \ref{dgnn_algo}. The other is the plane wave activation function, which associates with the DGTNN method introduced in section \ref{DGTNN}. Specially, a {\it discontinuous} plane wave neural network consists of a single hidden layer of $n\in \mathbb{R}$ neurons on each element $\Omega_k\in {\cal T}_h$. Define a function $\varphi^{\theta}: \mathbb{R}^d \rightarrow \mathbb{C}$ as follows:
\be \label{nn_app2_general}
\varphi^{\theta}({\bf x})|_{\Omega_k} = \sum_{j=1}^n c_j^{(k)} \text{e}^{\text{i} \omega {\bf d}_j^{(k)} \cdot {\bf x}}, ~~{\bf x}\in \Omega_k\quad (|{\bf d}_j^{(k)}| = 1),
\en
where $\{{\bf d}_j^{(k)}\in \mathbb{R}^d: |{\bf d}_j^{(k)}|=1\}$ are $n$ different propagation directions on
$k$-th element $\Omega_k\in {\cal T}_h$; $c_j^{(k)} \in \mathbb{C}$ are element-wise coefficients. Let $\text{T}_n({\mathcal T}_h)$ be the set of all functions of the form (\ref{nn_app2_general}), namely,
\be \label{vsign_general}
\text{T}_n({\mathcal T}_h) := \bigg\{v\in L^2(\Omega): ~v|_{\Omega_k}=\sum_{j=1}^n c_j^{(k)} \text{e}^{\text{i} \omega {\bf d}_j^{(k)} \cdot {\bf x}}, ~~{\bf x}\in \Omega_k,~~\Omega_k\in {\mathcal T}_h\bigg\}.
\en
%~~(c_j^{(k)} \in \mathbb{C}, ~{\bf d}_j^{(k)}\in \mathbb{R}^d, ~|{\bf d}_j^{(k)}| = 1)

The propagation directions $\{{\bf d}_j^{(k)}\}$ are initialized to be uniformly distributed on the unit circle.
The other parameters are set as follows.
\be \nonumber
 \omega=4\pi,
 \left\{\begin{array}{ll} h=\frac{1}{12}, n_r=2r+1, u_0=u^{(1)}_h & \text{for DGTNN}\\
 h=\frac{1}{20}, n_r=2r+9, u_0=0 & \text{for DGNN}.
\end{array}\right.
\en

\begin{figure}[H]
\begin{center}
\begin{tabular}{cc}
\epsfxsize=0.4\textwidth\epsffile{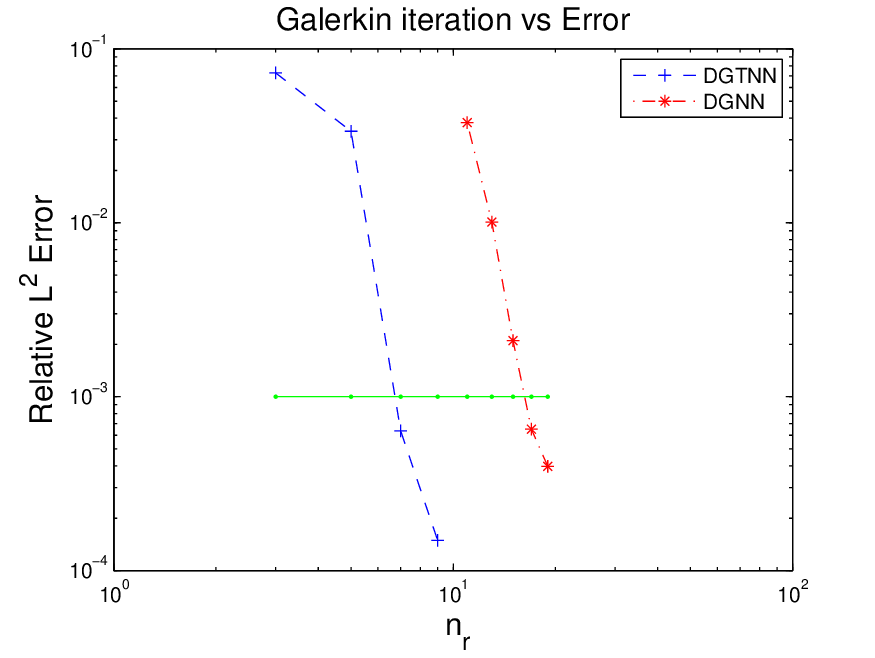}&
\epsfxsize=0.4\textwidth\epsffile{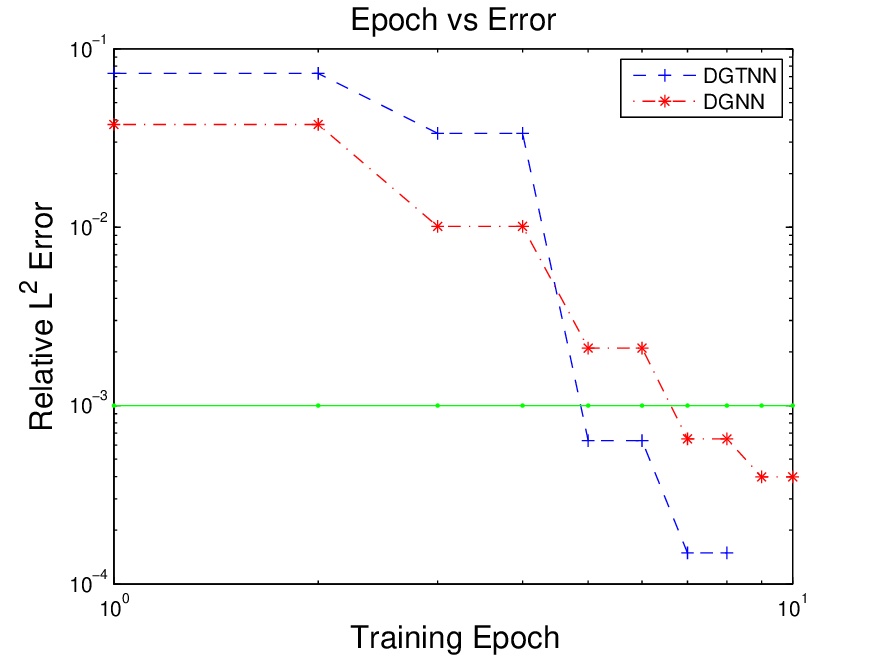}\\
\epsfxsize=0.4\textwidth\epsffile{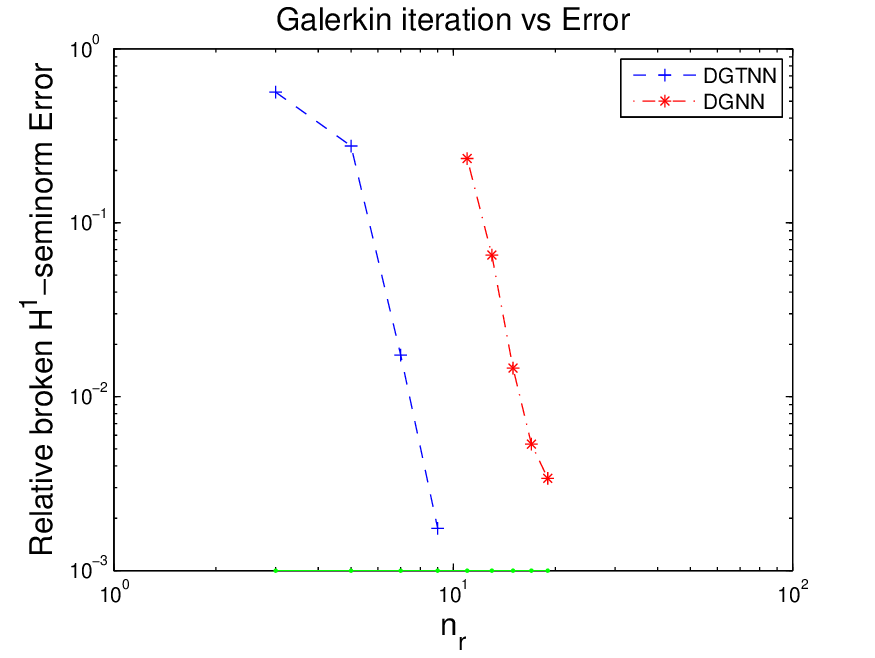}&
\epsfxsize=0.4\textwidth\epsffile{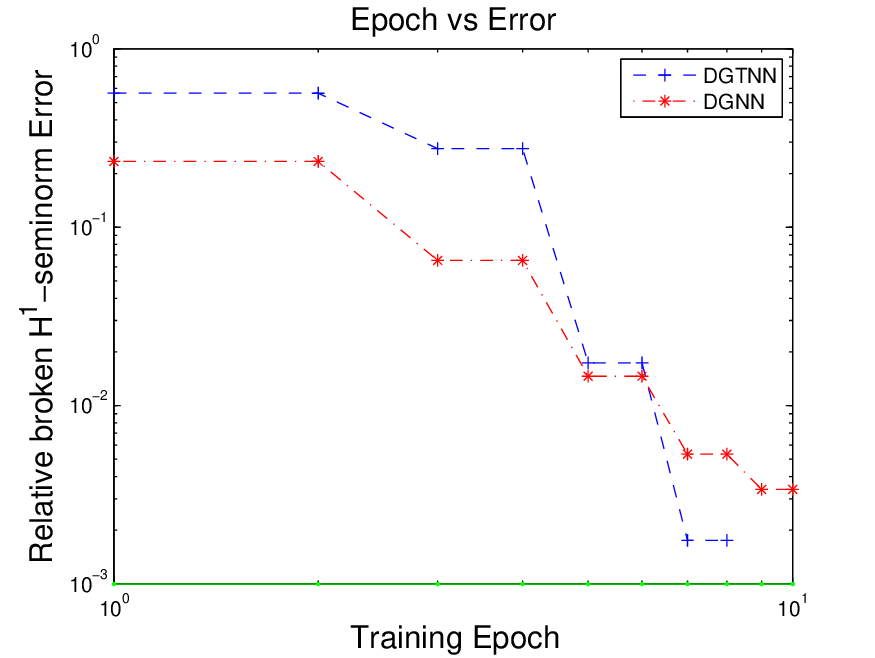}\\
\end{tabular}
\end{center}
 \caption{Displacement of a string (Helmholtz equation in two dimension). (Left-Up) Relative error in the $L^2$-norm at each Galerkin iteration. (Right-Up) The progress of the relative error in the $L^2$-norm within each Galerkin iteration. The x-axis thus denotes the cumulative training epoch over all Galerkin iterations. (Left-Bottom) Relative error in the broken $H^1$-seminorm at each Galerkin iteration. (Right-Bottom) The progress of the relative error in the broken $H^1$-seminorm within each Galerkin iteration. }
\label{2dnn1}
\end{figure}

Figure \ref{2dnn1} shows the relative errors of $u-u_r$ in the $L^2$-norm and the broken $H^1$-seminorm at each Galerkin iteration, respectively. We also provide the analogous results after each training epoch. Set the number of training epoches as 2 in Algorithm 4.1.

Compared with the single layer neural network architecture with general activation functions, Trefftz neural networks can achieve higher accuracy with smaller computational costs, that is, relatively coarser spatial width satisfying $\omega h=\frac{\pi}{3}$, and smaller number of neurons on each element.

\section{A high-accuracy DGNN method with two hidden layers for PDEs with variable coefficients}\label{twolayer_sec}

As we all know, Trefftz methods are viable only when the PDE is linear and its coefficients are piecewise-constant, which hinders the application of the DGTNN method to the models with variable coefficients, regardless of how many neurons are taken and how initial values are chosen. However, thanks to the flexibility and variability of neural network structures, we only need a general neural network structure only with two hidden layers to achieve a dramatic increase in efficiency and accuracy
when solving variable coefficient PDEs with the general source term.

Although the {\it generalized} Trefftz method has been used to solve {\it homogeneous} time-harmonic Helmholtz equations \cite{imbertnume,imbertscal}, time-domain acoustic equations \cite{imbertmosto}, %Schr\"{o}dinger equations,
 and time-harmonic electromagnetic field equations \cite{imbertmax} with variable coefficients, its construction of generalized Trefftz basis functions maps complex nonlinear relationships. Therefore, it is not wise to introduce neural networks to train their nonlinear parameters. Otherwise it incurs significant computational costs. Instead, we propose a simple, economical, and high-precision {\it two-hidden-layer} DGNN algorithm.  %\textcolor{red}{For convenience, we call the method as DGNN-TWO method.}

 The two-hidden-layer feed-forward discontinuous neural network is defined as
 \be
 u_{NN}({\bm x},t; {\bm W}, {\bm b})|_{Q_k} = \mathcal{L}^{(k)} \circ T_2^{(k)} \circ T_1^{(k)}({\bm x},t).
 \en
In each hidden layer $i=1,2$, the nonlinear mapping is ${\bf Z}^{(k)}_i=T_i^{(k)}({\bf U}^{(k)}_i)=\sigma({\bf U}^{(k)}_i)=\sigma({\bf W}^{(k)}_i\times {\bf Z}^{(k)}_{i-1}+{\bf b}^{(k)}_i)$ with weights ${\bf W}^{(k)}_i\in \mathbb{C}^{m_i\times m_{i-1}}$ and biases ${\bf b}^{(k)}_i\in \mathbb{C}^{m_i}$, where ${\bf Z}^{(k)}_0=({\bm x},t)^t|_{Q_k}$ (The superscript ``t" denote the transpose of a vector), $m_0=d+1$ is the input dimension. In the output layer, $\mathcal{L}^{(k)}: \mathbb{C}^{m_2}\rightarrow \mathbb{C}$ is the linear mapping, namely, $u_{NN}({\bm x},t)|_{Q_k} = {\bf W}^{(k)}_3 \times {\bf Z}^{(k)}_2$.

\begin{figure}[H]
%\vspace{-2cm}thb
\begin{center}
\begin{tabular}{c}
\epsfxsize=0.8\textwidth\epsffile{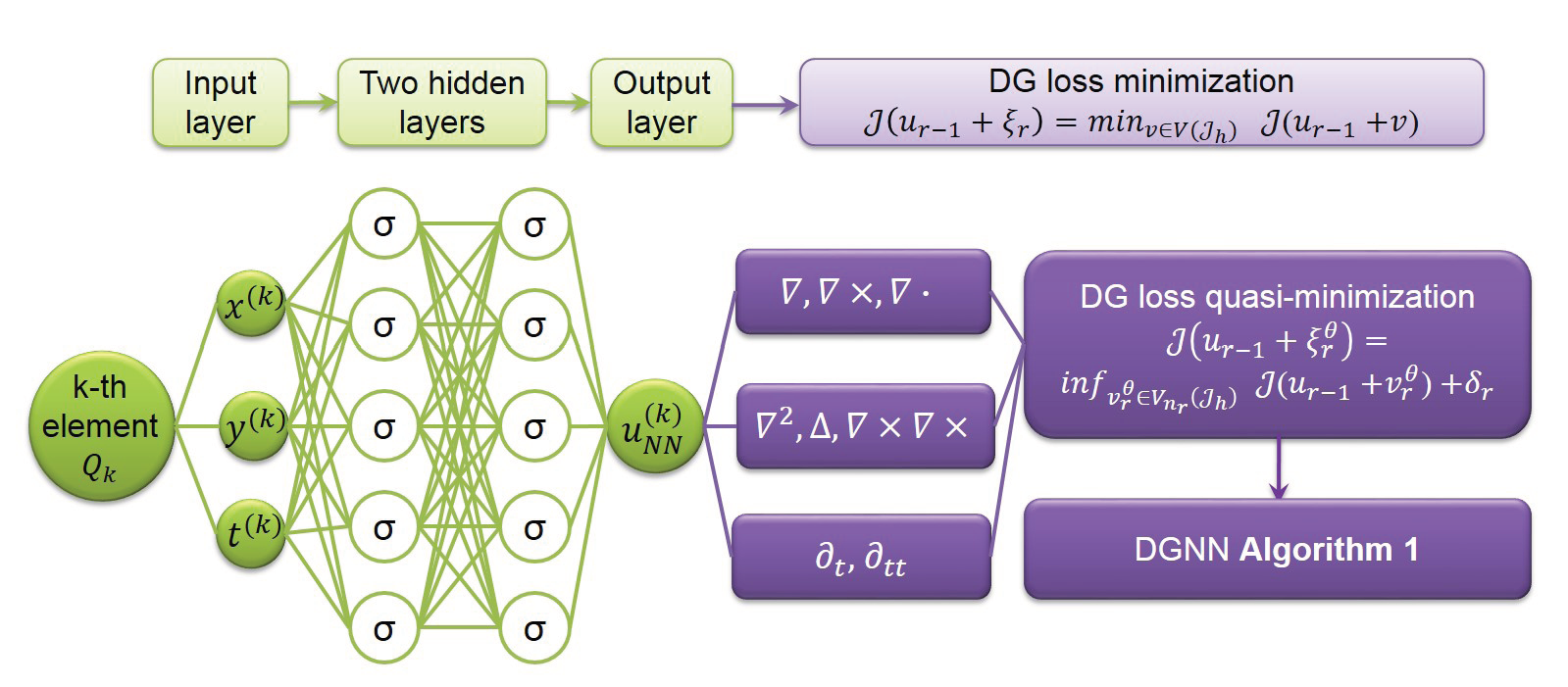}\\
\end{tabular}
\end{center}
 \caption{ A schematic for the discontinuous Galerkin Neural network method with two hidden layers. }
\label{2hlnn}
\end{figure}

Figure \ref{2hlnn} shows a schematic for our DGNN method with two hidden layers. % in a DGNN iterative method.

\subsection{Helmholtz equations with variable coefficient $\rho = x^2+y^2$}\label{variahelm}

Consider the Helmholtz equation (\ref{helmmodel}) with a variable coefficient $\rho = x^2+y^2$.
The exact solution of the above problem is defined in the closed form (see \cite{yhijcom})
$$u_{ex}(x,y) = y~\text{sin}(\omega x)$$
with $\Omega=[0,1]\times[0,1]$, and $g=({\partial \over
\partial {\bf n}} + \text{i}\omega)u_{ex}$.
Then the source term $f=\omega^2 y ~\text{sin}(\omega x) (1-\rho).$ The operator $\mathcal{A}$ is defined by $\mathcal{A}= -\Delta -\omega^2\rho$, and the other operators $\mathcal{B}, \mathcal{C}_1, \mathcal{C}_2$ are defined by (\ref{foperadefin}), respectively.

Set $u_0=0$, $\omega=4\pi, \omega h=\pi/3$, $m_1=7$, and $m_2 = n_r=2r+9$ for each Galerkin iteration $r$.

Initialize nonlinear parameters as follows. Weights ${\bf W}^{(k)}_1\in \mathbb{C}^{m_1\times 2}$ in the first hidden layer are uniformly distributed on the unit circle.
\be\nonumber
{\bf W}^{(k)}_1(j,:)=(cos\theta_j, sin\theta_j), ~~~\theta_j=\frac{2\pi}{m_1}j, j=1,\cdots,m_1.
\en
Set biases ${\bf b}^{(k)}_1={\bf 1}\in \mathbb{C}^{m_1}$. Weights ${\bf W}^{(k)}_2\in \mathbb{C}^{m_2\times m_1}$ in the second hidden layer are set as
${\bf W}^{(k)}_2 = \text{i}\omega Q$, where $m_2\geq m_1$ (Note that $\text{i}$ is the pure imaginary unit), and the matrix $Q$ is composed of the first $m_1$ columns of a $m_2$-order random orthogonal matrix. Set biases ${\bf b}^{(k)}_2={\bf 1}\in \mathbb{C}^{m_2}$.

We want to explain why the nonlinear parameters ${\bf W}^{(k)}_i$ ($i=1,2$) are initialized in this model-driven manner. Without causing confusion, we will omit the superscript $k$.
Set ${\bf Z}_2 = (Z_{2,\ell})_{\ell=1}^{m_2}$, similarly for ${\bf U}_2, {\bf U}_1$. Set ${\bf W}_2 = (W^{(2)}_{\ell s})$, $1\leq \ell \leq m_2, 1 \leq s \leq m_1$, similarly for ${\bf W}_1 = (W^{(1)}_{st})$.

By the chain rule, we have
\beq \label{laprole}
\Delta {Z}_{2,\ell } & = \sigma''({U}_{2,\ell }) \sum_{t=1}^d \bigg( \sum_{s=1}^{m_1} W^{(2)}_{\ell s}\sigma'(U_{1,s})W^{(1)}_{st} \bigg)^2
\cr
& + \sigma'({U}_{2,\ell })  \sum_{s=1}^{m_1} W^{(2)}_{\ell s}\sigma''(U_{1,s}) \sum_{t=1}^d \big( W^{(1)}_{st} \big)^2.
\eq

For the adopted sigmoid activation function $\sigma(\xi)=\frac{1}{1+e^{-\xi}}, \xi\in \mathbb{C}, \xi\neq \text{i}(\pi+2\pi \upsilon), \upsilon\in \mathbb{N}$, one may be quick to calculate $\sigma'(\xi)=\sigma(1-\sigma)$ and $\sigma''(\xi)=\sigma(1-\sigma)(1-2\sigma)$. It can be verified that selecting {\it pure imaginary matrix} ${\bf W}^{(k)}_2$ and ${\bf b}^{(k)}_2={\bf 1}$ ensures the boundedness of the functions $\sigma'({U}_{2,\ell })$ and $\sigma''({U}_{2,\ell })$. Of course, $\sigma'(U_{1,s})$ and $\sigma''(U_{1,s})$ are bounded owing to $U_{1,s}\in \mathbb{R}$.
By the Cauchy-Schwarz inequality, we can bound the first term of the sum (\ref{laprole}) by
\be
\bigg|\sigma''({U}_{2,\ell }) \sum_{t=1}^d \bigg( \sum_{s=1}^{m_1} W^{(2)}_{\ell s}\sigma'(U_{1,s})W^{(1)}_{st} \bigg)^2 \bigg| \leq C \bigg( \sum_{s=1}^{m_1} |W^{(2)}_{\ell s}|^2 \bigg)
\bigg( \sum_{t=1}^d \sum_{s=1}^{m_1} |W^{(1)}_{st}|^2 \bigg),
\en
where $C$ denotes a constant independent of any variable. Obviously, the second term of the sum (\ref{laprole}) can be bounded by
\be
\bigg|\sigma'({U}_{2,\ell})  \sum_{s=1}^{m_1} W^{(2)}_{\ell s}\sigma''(U_{1,s}) \sum_{t=1}^d \big( W^{(1)}_{st} \big)^2 \bigg| \leq C \bigg( \sum_{s=1}^{m_1} |W^{(2)}_{\ell s}| \bigg)
 \bigg( \sum_{t=1}^d \sum_{s=1}^{m_1} |W^{(1)}_{st}|^2 \bigg).
\en

Comparing the two inequalities above, we can deduce that the first term of of the sum of two terms in (\ref{laprole}) is dominant in the Laplace operator on ${Z}_{2,\ell }$, and it follows the order of $\mathcal{O}(\omega^2)$ by the original Helmholtz equation $-\Delta u-\omega^2 \rho u=f$. Meanwhile, the second term of the sum occupies a recessive role, which is about the order of $\mathcal{O}(\omega)$. The above findings determine the choice of initial values of nonlinear parameters: the row vectors of the matrix ${\bf W}_1$ are unit vectors, and ${\bf W}_2= \text{i} \mathcal{O}(\omega)$.

\begin{figure}[H]
%\vspace{-2cm}thb
\begin{center}
\begin{tabular}{cc}
\epsfxsize=0.4\textwidth\epsffile{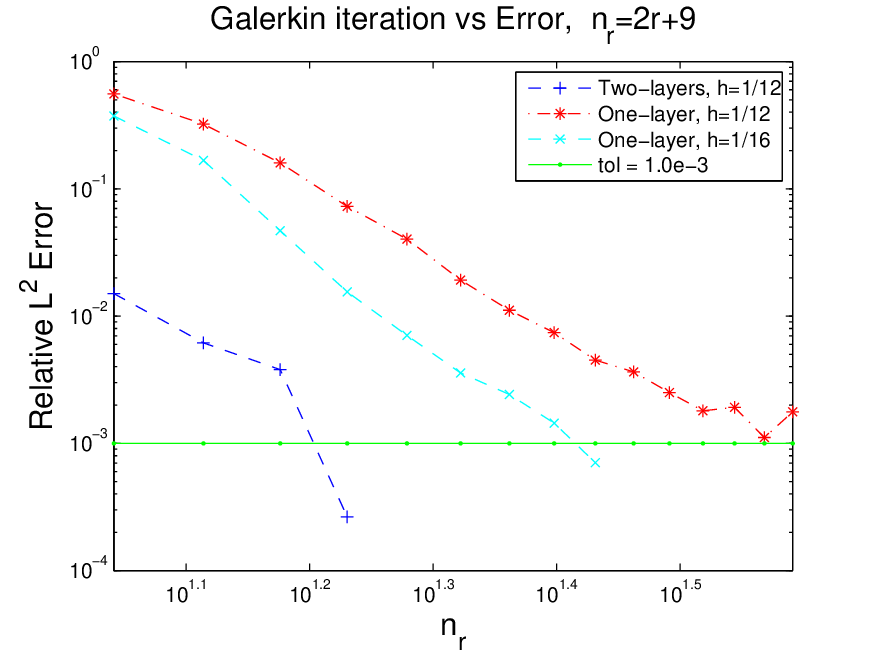}&
\epsfxsize=0.4\textwidth\epsffile{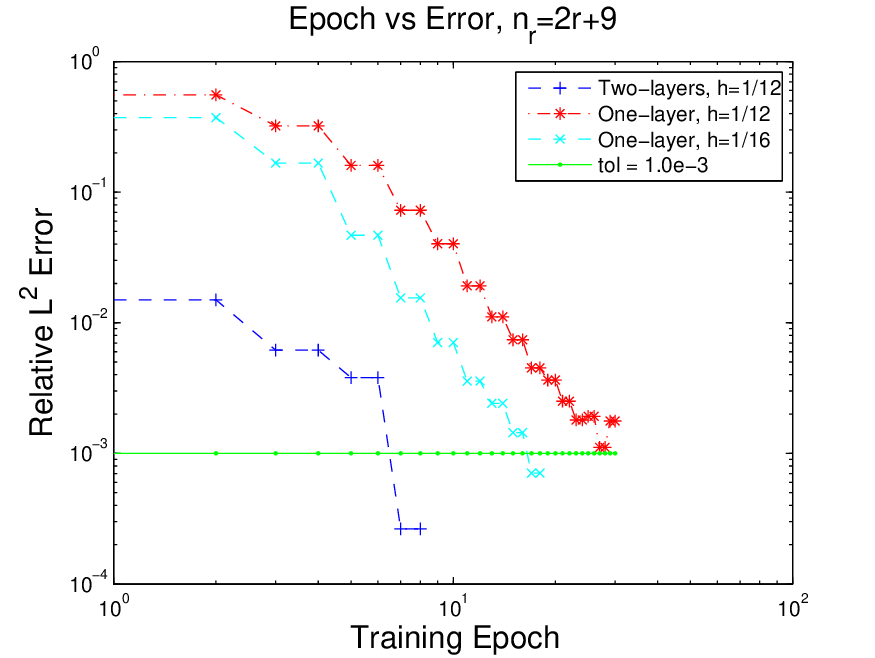}\\
\epsfxsize=0.4\textwidth\epsffile{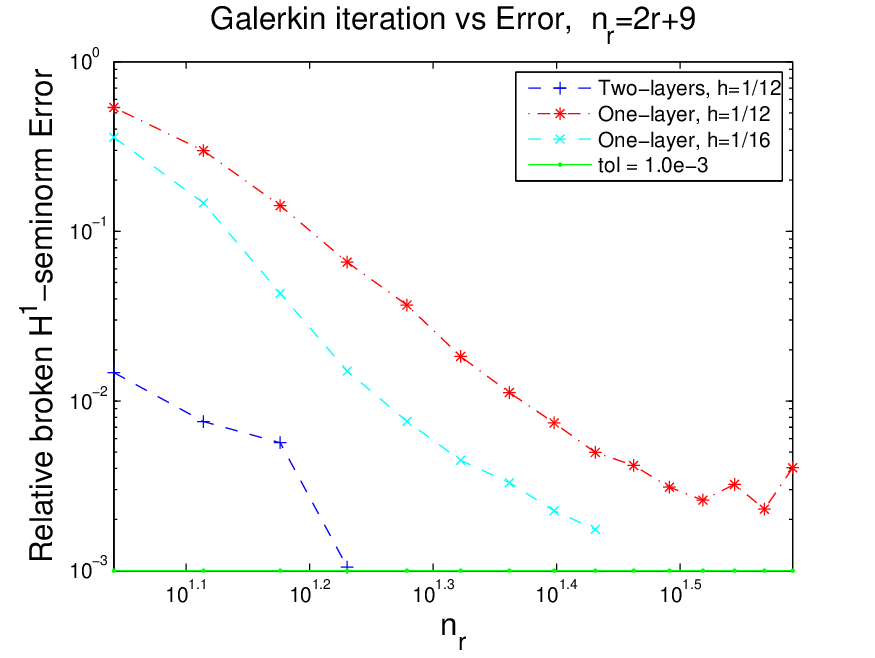}&
\epsfxsize=0.4\textwidth\epsffile{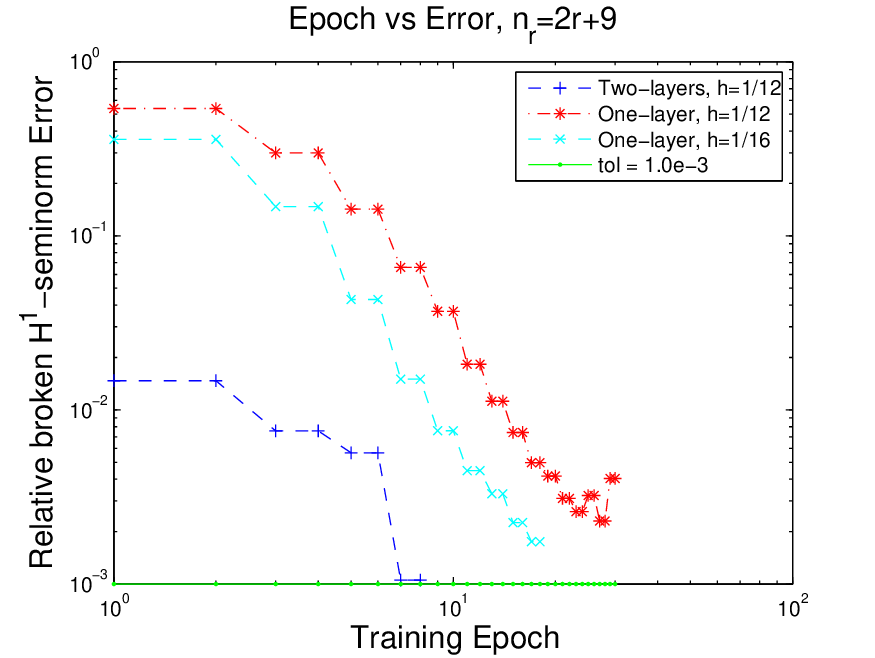}\\
\end{tabular}
\end{center}
 \caption{Displacement of a string (Helmholtz equation with variable coefficients). (Left-Up) Relative error in the $L^2$-norm at each Galerkin iteration. (Right-Up) The progress of the relative error in the $L^2$-norm within each Galerkin iteration. The x-axis thus denotes the cumulative training epoch over all Galerkin iterations. (Left-Bottom) Relative error in the broken $H^1$-seminorm at each Galerkin iteration. (Right-Bottom) The progress of the relative error in the broken $H^1$-seminorm within each Galerkin iteration. }
\label{2dvarihelm_nn1}
\end{figure}

%\textcolor{red}{
Table \ref{2dconsthelmtable} shows computational complexity of the DGNN methods with one-hidden layer and two-hidden layers, respectively, for different $h$ and neurons. Here ``DOFs" represents the number of degree of freedoms of the least square linear system (\ref{dresivari}). Figure \ref{2dvarihelm_nn1} shows the relative errors of $u-u_r$ in the $L^2$-norm and the broken $H^1$-seminorm at each Galerkin iteration, respectively. We also provide the analogous results after each training epoch. Set the number of training epoches as 2 in Algorithm 4.1.

\begin{center}
       \tabcaption{}%\vskip -0.3in
\label{2dconsthelmtable}
      Comparisons of computational complexity.  \vskip 0.1in
\begin{tabular}{|c|c|c|c|} \hline
   \text{Method} & \text{Neurons} & \text{DOFs} & \text{Relative} $L^2$ \text{Error} \\ \hline
  \text{Two-layers case}, $\omega h=\frac{\pi}{3}$ & 3456 & 2448 & 2.64e-4  \\ \hline
  \text{One-layer case}, $\omega h=\frac{\pi}{3}$ & 5616 & 5616 & 1.77e-3 \\  \hline
  \text{One-layer case}, $\omega h=\frac{\pi}{4}$ & 6912 & 6912 &  7.05e-4 \\  \hline
   \end{tabular}
     \end{center}

Numerical results validate that, the resulting approximate solutions generated by the DGNN method employing two-hidden-layer neural network architectures can reach the highest accuracy. Moreover, the cost of computing linear coefficients ${\bf W}^{(k)}_3$ by Galerkin least squares (\ref{dresivari}) and updating nonlinear weights ${\bf W}^{(k)}_1, {\bf W}^{(k)}_2$ by the  gradient descent algorithm (\ref{gradientIMPLE}) is the smallest.

\begin{remark}
It is gratifying that, under the selection of the default initial value $u_0=0$, thanks to the reasonable selection of nonlinear initial values ${\bf W}$ and {\bf b},
the DGNN algorithm can ensure the fast convergence, whether the discontinuous neural network activated by general smooth functions such as sigmoid functions is a single hidden layer or a double hidden layer structure.
 It should be emphasized that the template for selecting the initial values of nonlinear weights based on the considered equations determines that the DGNN algorithm is driven by the model.
\end{remark}

\section{Numerical experiments} \label{numerical}

In the remaining benchmark experiments, we will test the effectiveness and accuracy of the proposed DGNN method for the time-harmonic and time-dependent wave equations. We would like to point out that, the DGNN method reduces the relative $L^2$ error by at least one order of magnitude under very economical conditions, compared to the existing PINN method and its variations.

Set $tol = 10^{-3}$, the number of training epoches as 2 in Algorithm 4.1.

\subsection{Time-harmonic Maxwell's equations} \label{numesec}
Consider the Maxwell equations with the lowest order absorbing boundary condition in three space dimensions \cite{ALZOU,BZZOU}:
\begin{equation} \label{nnmaxeq}
\left\{ \begin{aligned}
     &  \nabla\times(\frac{1}{i\omega\mu}\nabla\times {\bf E}) + \text{i}\omega\varepsilon {\bf E} = {\bf f} & \text{in}\quad
\Omega,\\
    & -{\bf E}\times {\bf n}+\frac{\sigma}{\text{i}\omega\mu}((\nabla\times {\bf E})\times  {\bf n})\times {\bf n}={\bf g} & \text{on}\quad\Gamma.
                          \end{aligned} \right.
                          \end{equation}
Here $\omega>0$ is the temporal frequency of the field, ${\bf f} \in {\bf L}^2(\Omega)^3$ and ${\bf g}\in {\bf L}_{\bf T}^2(\partial\Omega)^3$. The material coefficients
$\varepsilon,\mu$ and $\sigma$ are assumed to be piecewise constant in the
whole domain. In particular, if $\varepsilon$ takes complex valued,
then the material is an absorbing medium; else the material is a non-absorbing medium (see \cite{hmm}).

The operators $\mathcal{A}, \mathcal{B}, \mathcal{C}_1, \mathcal{C}_2$ are defined by, respectively,
\be
\mathcal{A}{\bf E } = \nabla\times(\frac{1}{\text{i}\omega\mu}\nabla\times {\bf E}) + \text{i}\omega\varepsilon {\bf E}, ~\mathcal{B} {\bf E } = -{\bf E}\times {\bf n}+\frac{\sigma}{\text{i}\omega\mu}((\nabla\times {\bf E})\times  {\bf n})\times {\bf n},
\en
\be
 \mathcal{C}_1 ({\bf E })  = \llbracket{ {\bf E}} \rrbracket_{\bf T} ~\text{and}
 ~ \mathcal{C}_2 ({\bf E }) = \llbracket {1\over \text{i}\omega\mu}(\nabla\times{{\bf E}}) \rrbracket_{\bf T} ~\text{over}~
\mathcal{F}_h^{\text{I}}.
\en

The relaxation parameters $\{\lambda_i\}$ are initialized as follows.
\be\nonumber
 \lambda_0=1, ~\lambda_1=\lambda_2=\lambda_3=\omega^2.
\en

\subsubsection{The case of constant coefficients}
We test our method with the material coefficients $\mu =1, \varepsilon = 1 + \text{i}$ and $\sigma=\sqrt{\frac{\mu}{|\varepsilon|}}$ according to the known
formula introduced in \cite[section 1]{hmm}. Consider the following analytical solution
\begin{equation}
\label{eq43} {\bf E}_{\text{ex}}=
\varepsilon\omega(xz~\text{cos}y,-z~\text{sin}y, xy)^t
\end{equation}
in a space domain $\Omega=(0,1)\times(0,1)\times(0,1)$.

We employ both the DGNN method introduced in section \ref{helmdis} and the DGTNN method introduced in section \ref{constmodel} to solve
(\ref{nnmaxeq}). Two types of {\it vector} activation functions are adopted to generate neural network spaces. One is the constant {\it vector} multiplied by the sigmoid activation functions. For the case, Initialize weights as ${\bf W}^{(k)} = \text{i}\kappa \tilde{\bf W}^{(k)}\in \mathbb{C}^{n\times 3}$, where $\tilde{\bf W}^{(k)}$ are initialized by the spherical codes from \cite{refsite}. Set biases ${\bf b}^{(k)}={\bf 1}\in \mathbb{C}^{n}$.

The other is the polarization {\it vector} multiplied by the plane wave activation functions, which associates with the proposed DGTNN method introduced in section \ref{DGTNN}. Specially, a {\it discontinuous} plane wave neural network consists of a single hidden layer of $2n\in \mathbb{R}$ neurons on each element $\Omega_k\in {\cal T}_h$. Define a function $\varphi^{\theta}: \mathbb{R}^d \rightarrow \mathbb{C}^d$ as follows:
\be \label{nn_app2_generalmax}
\varphi^{\theta}({\bf x})|_{\Omega_k} = \sum_{j=1}^n \big( c_j^{(k)} \sqrt{\mu}~{\bf f}_{j}^{(k)}~e^{\text{i}\kappa {\bf d}_{j}^{(k)}\cdot{\bf x}}
 + c_{j+n}^{(k)} \sqrt{\mu}~ {\bf f}_{j+n}^{(k)}~e^{\text{i}\kappa {\bf d}_{j}^{(k)}\cdot{\bf x}} \big), ~~{\bf x}\in \Omega_k,\
%c_j^{(k)} \text{e}^{\text{i} \omega {\bf d}_j^{(k)} \cdot {\bf x}}, ~~{\bf x}\in \Omega_k\quad (|{\bf d}_j^{(k)}| = 1),
\en
where $c_j^{(k)} \in \mathbb{C}$ are element-wise coefficients, $\kappa=\omega\sqrt{\mu\varepsilon}$, $\{{\bf d}_j^{(k)}\in \mathbb{R}^d: |{\bf d}_j^{(k)}|=1\}$ are $n$ different propagation directions on $k$-th element $\Omega_k\in {\cal T}_h$. For $j=1,\cdots,n$, every unit real polarization vector ${\bf f}_j^{(k)}$ is orthogonal to ${\bf d}_j^{(k)}$, and ${\bf f}_{j+n}^{(k)}={\bf f}_j^{(k)}\times {\bf d}_j^{(k)}$. Let ${\bf T}_{2n}({\mathcal T}_h)$ be the set of all functions of the form (\ref{nn_app2_generalmax}), namely,
\be \label{vsign_generalmax}
{\bf T}_{2n}({\mathcal T}_h) := \bigg\{ {\bf F} \in {\bf L}^2(\Omega): ~{\bf F}|_{\Omega_k}=\sum_{j=1}^n \big( c_j^{(k)} \sqrt{\mu}~{\bf f}_{j}^{(k)}~e^{\text{i}\kappa {\bf d}_{j}^{(k)}\cdot{\bf x}}
 + c_{j+n}^{(k)} \sqrt{\mu}~ {\bf f}_{j+n}^{(k)}~e^{\text{i}\kappa {\bf d}_{j}^{(k)}\cdot{\bf x}} \big) \bigg\}.
\en
The propagation directions $\{{\bf d}_j^{(k)}\}$ are initialized by the spherical codes from \cite{refsite}.

For the DGTNN method, as introduced in the previous section \ref{local_spec} for computing the initial value ${\bf E}_0$, the variational problem of \((\ref{cha4.n1})\) is to find \(
{\bf E}^{(1)}_k\in H(\text{curl}, \Omega^{\ast}_k)\) such that
\begin{equation} \label{eq6.3}
\left\{ \begin{aligned}
     &\int_{\Omega^{\ast}_k}(\frac{1}{\text{i}\omega\mu}\nabla\times{\bf E}^{(1)}_k\cdot \overline{\nabla\times{\bf
F}_k}   + \text{i}\omega\varepsilon{\bf E}^{(1)}_k \cdot \overline{{\bf F}_k} )dV+
\int_{\partial\Omega^{\ast}_k}\frac{1}{\sigma}({\bf E}^{(1)}_k\times{\bf
n})\times{\bf n}\cdot \overline{{\bf F}_k} dS=
\int_{\Omega^{\ast}_k}{\bf f} \cdot \overline{{\bf F}_k} dV,\\
    & \quad\quad\quad\quad\quad\quad\quad\quad\quad\forall{\bf F}_k\in
{\bf \mbox{H}}(\text{curl}, \Omega^{\ast}_k)~(k = 1,2,\ldots,N).
                          \end{aligned} \right.
                          \end{equation}

We choose local spectral elements of order $m=3$ as the initial value ${\bf E}_0$ in the original DGTNN iterative algorithm. The discrete variational problems of (\ref{eq6.3}) are: to find ${\bf E}^{(1)}_{k,h}\in {\bf S}_3(\Omega^{\ast}_k)$ such that
\begin{equation} \label{3.new2}
\left\{ \begin{aligned}
     &\int_{\Omega^{\ast}_k}(\frac{1}{\text{i}\omega\mu}\nabla\times{\bf E}^{(1)}_{k,h} \cdot \overline{\nabla\times{\bf
F}_k}  + \text{i}\omega\varepsilon{\bf E}^{(1)}_{k,h} \cdot \overline{{\bf F}_k}) dV +
\int_{\partial\Omega^{\ast}_k}\frac{1}{\sigma}({\bf E}^{(1)}_{k,h} \times{\bf
n})\times{\bf n} \cdot \overline{{\bf F}_k} dS =
\int_{\Omega^{\ast}_k}{\bf f} \cdot \overline{{\bf F}_k} dV ,\\
    & \quad\quad\quad\quad\quad\quad\quad\quad\quad\forall{\bf F}_k\in
{\bf S}_3(\Omega^{\ast}_k)~(k = 1,2,\ldots,N).
                          \end{aligned} \right.
                          \end{equation}

The other parameters are set as follows.
\be \nonumber
 \omega=4\pi,
 \left\{\begin{array}{ll} h=\frac{1}{12}, n_r=(r+1)^2, {\bf E}_0= {\bf E}^{(1)}_h & \text{for DGTNN}\\
h=\frac{1}{20}, n_r=(r+1)^2, {\bf E}_0={\bf 0} & \text{for DGNN}.
\end{array}\right.
\en

\begin{figure}[H]
%\vspace{-2cm}thb
\begin{center}
\begin{tabular}{cc}
\epsfxsize=0.4\textwidth\epsffile{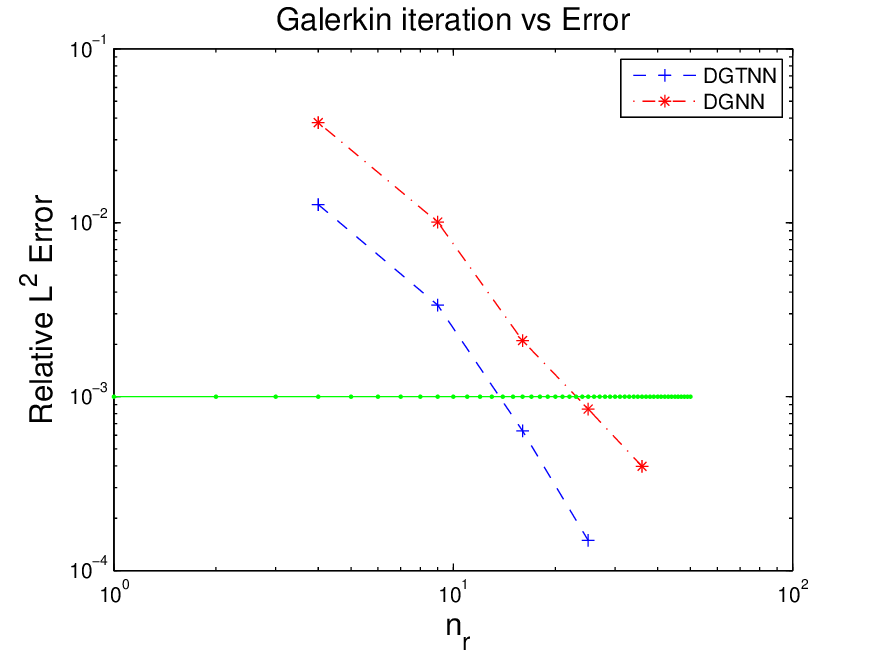}&
%COMPARE_2dhelm_galer_l2.eps
\epsfxsize=0.4\textwidth\epsffile{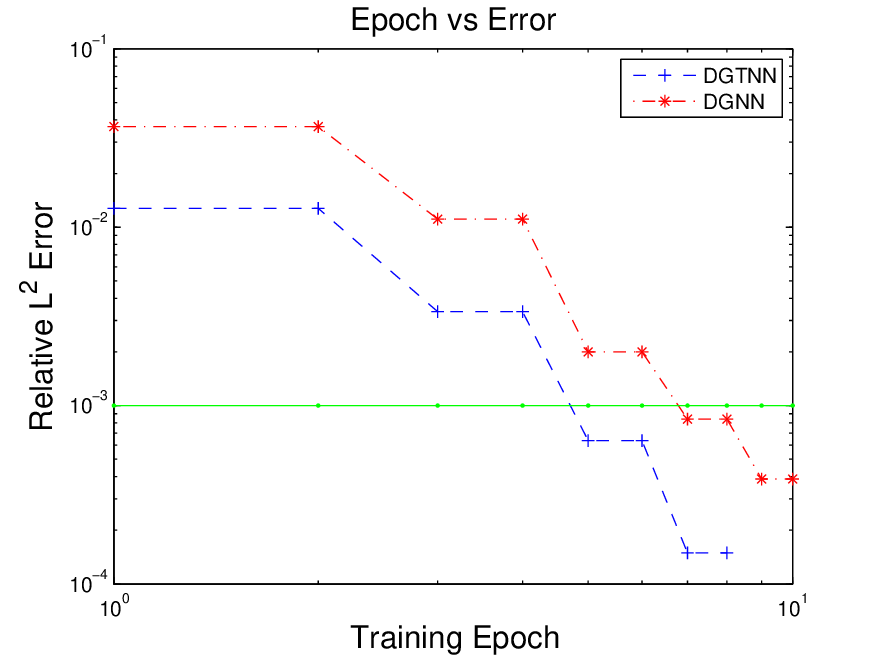}\\
\epsfxsize=0.4\textwidth\epsffile{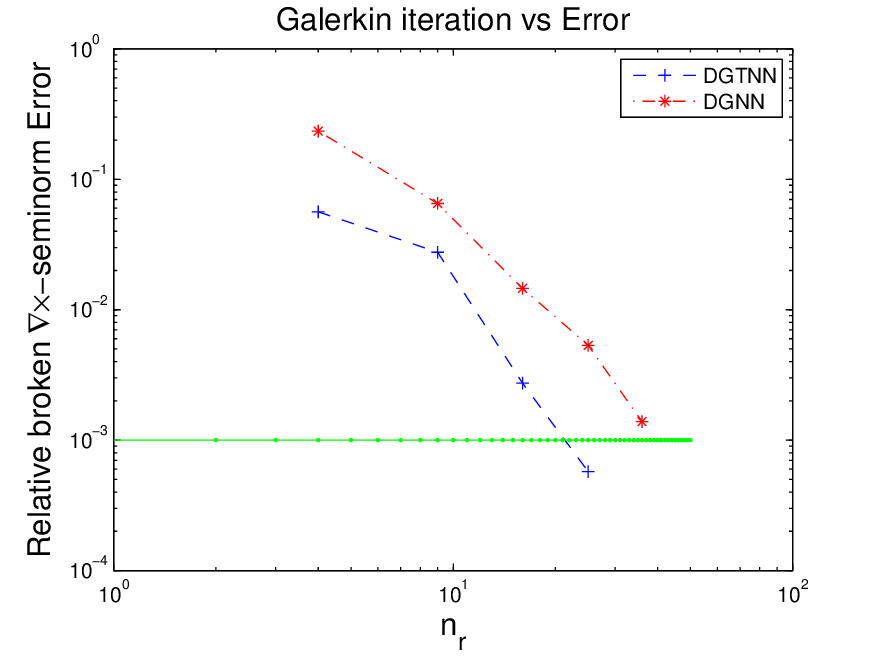}&
\epsfxsize=0.4\textwidth\epsffile{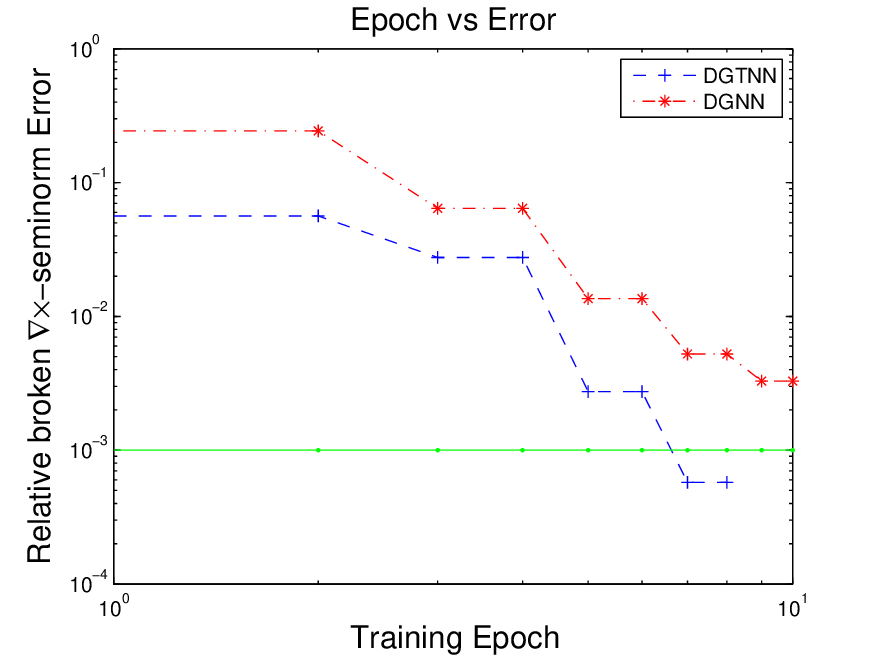}\\
\end{tabular}
\end{center}
 \caption{Displacement of a string (Maxwell's equations in three dimension). (Left-Up) Relative error in the $L^2$-norm at each Galerkin iteration. (Right-Up) The progress of the relative error in the $L^2$-norm within each Galerkin iteration. The x-axis thus denotes the cumulative training epoch over all Galerkin iterations. (Left-Bottom) Relative error in the broken $\nabla\times$-seminorm at each Galerkin iteration. (Right-Bottom) The progress of the relative error in the broken $\nabla\times$-seminorm within each Galerkin iteration. }
\label{3dconsmaxnn1}
\end{figure}

Figure \ref{3dconsmaxnn1} shows the relative errors of ${\bf E}-{\bf E}_r$ in the $L^2$-norm and the broken $\nabla\times$-seminorm at each Galerkin iteration, respectively. We also provide the analogous results after each training epoch. %Set the number of training epoches as 2 in {\bf Algorithm 4.1}.

%\textcolor{red}{insert results of $u_0=0$}

Compared with the single layer neural network architecture with general activation functions, Trefftz neural networks can achieve higher accuracy with smaller computational costs, that is, relatively coarser spatial width satisfying $\omega h=\frac{\pi}{3}$, and smaller number of neurons on each element.

\subsubsection{The case of variable coefficients}

We still consider the analytical solution (\ref{eq43}), but with the different material coefficients $\varepsilon = 1 + (x+y+z)\text{i}$ and
$\sigma=\sqrt{\frac{\mu}{|\varepsilon|}}$.

Initialize nonlinear parameters as follows. Weights ${\bf W}^{(k)}_1\in \mathbb{C}^{m_1\times 3}$ in the first hidden layer are generated by the spherical codes from \cite{refsite}. Set biases ${\bf b}^{(k)}_1={\bf 1}\in \mathbb{C}^{m_1}$. Weights ${\bf W}^{(k)}_2\in \mathbb{C}^{m_2\times m_1}$ in the second hidden layer are set as
${\bf W}^{(k)}_2 = \text{i}\kappa Q$, where $m_2\geq m_1$, the matrix $Q$ is composed of the first $m_1$ columns of a $m_2$-order random orthogonal matrix. Set biases ${\bf b}^{(k)}_2={\bf 1}\in \mathbb{C}^{m_2}$. The reason why the nonlinear parameters ${\bf W}^{(k)}_i$ ($i=1,2$) are initialized in this model-driven manner can refer to section \ref{variahelm}.

Set ${\bf E}_0=0$, $\omega=4\pi, \omega h=\pi/3$, $m_1=16$, and $m_2 = n_r = m_1+4(r-1)$ for each Galerkin iteration $r$.

\begin{figure}[H]
%\vspace{-2cm}thb
\begin{center}
\begin{tabular}{cc}
\epsfxsize=0.4\textwidth\epsffile{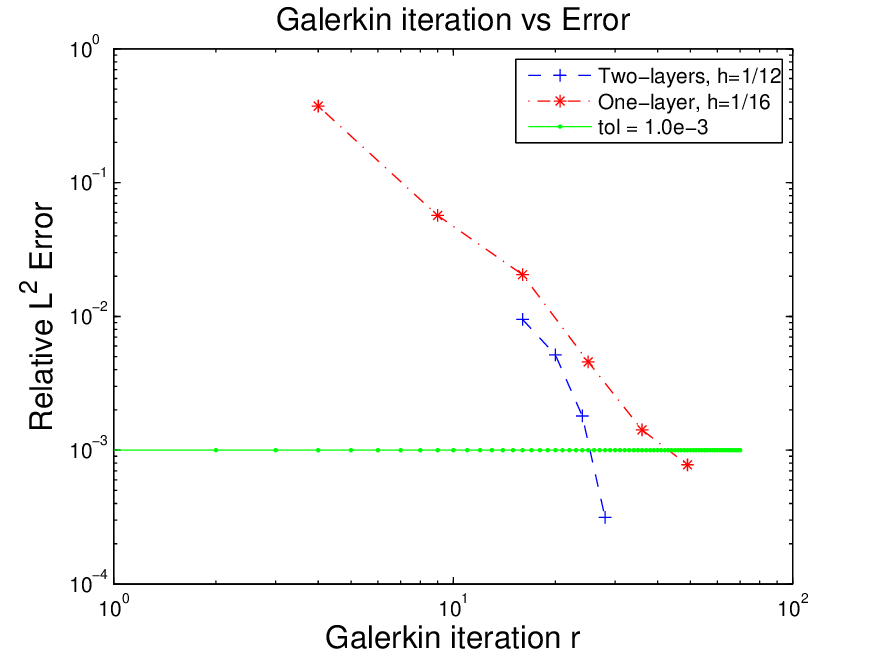}&
\epsfxsize=0.4\textwidth\epsffile{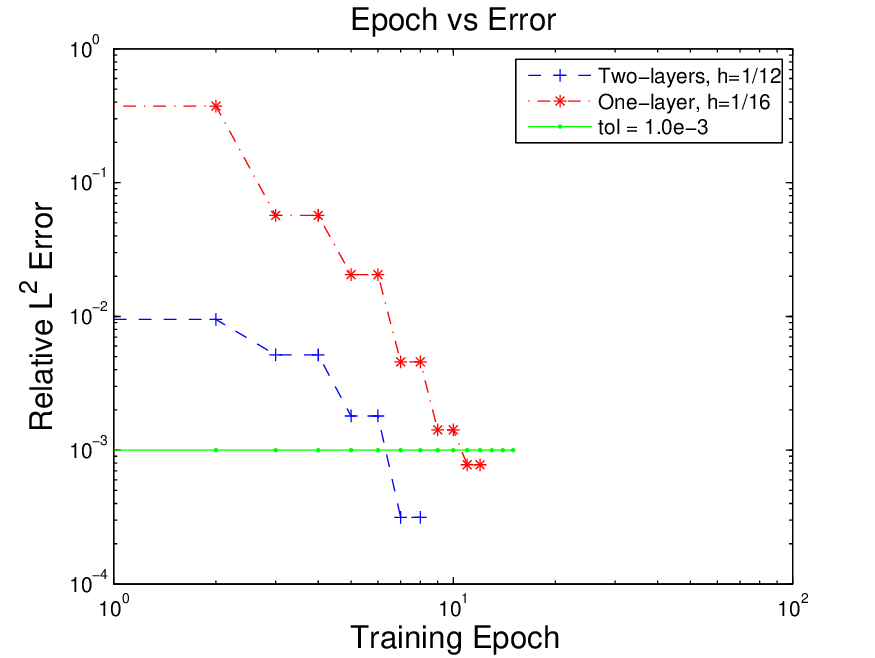}\\
\epsfxsize=0.4\textwidth\epsffile{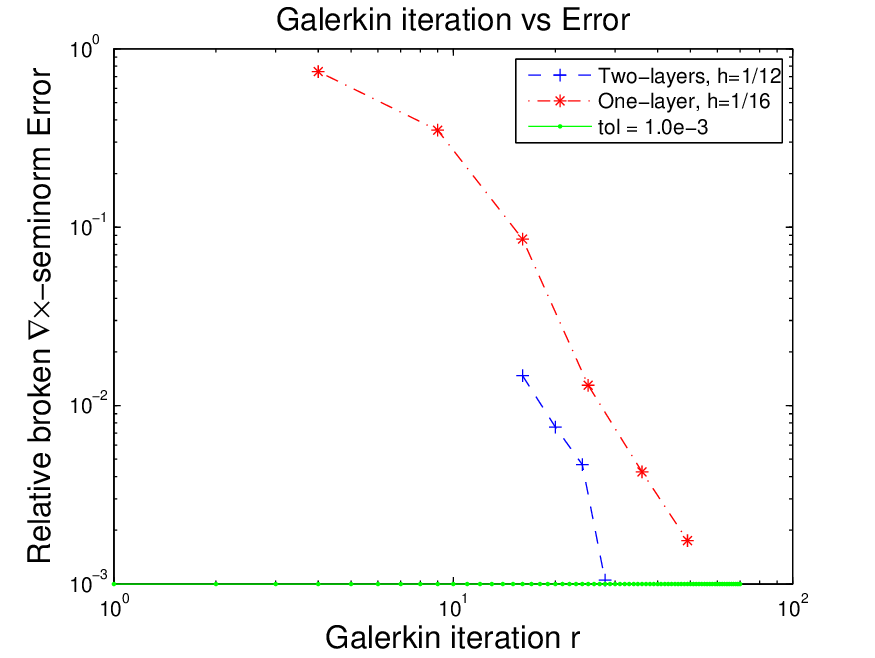}&
\epsfxsize=0.4\textwidth\epsffile{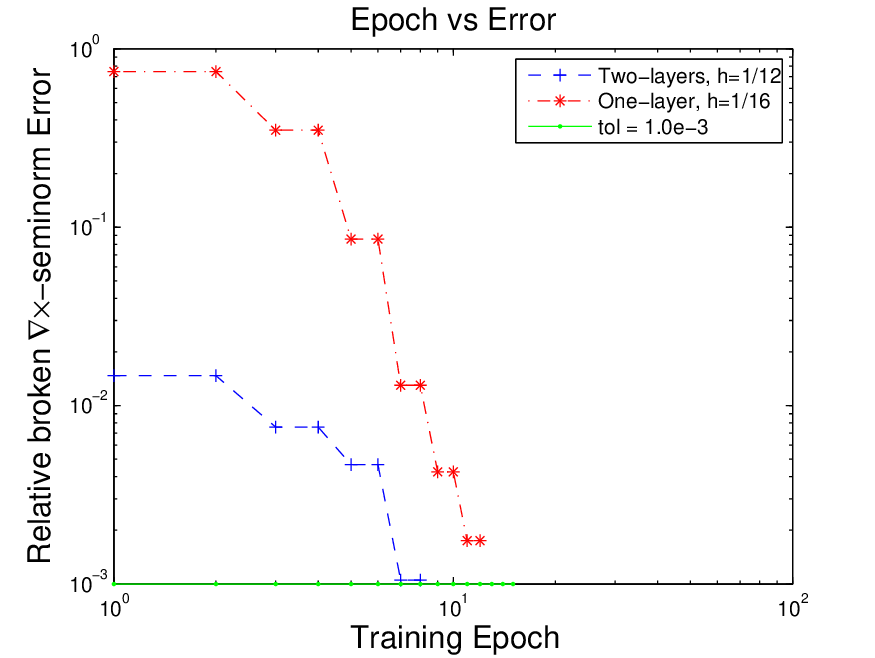}\\
\end{tabular}
\end{center}
 \caption{Displacement of a string (Maxwell's equations in three dimension). (Left-Up) Relative error in the $L^2$-norm at each Galerkin iteration. (Right-Up) The progress of the relative error in the $L^2$-norm within each Galerkin iteration. The x-axis thus denotes the cumulative training epoch over all Galerkin iterations. (Left-Bottom) Relative error in the broken $\nabla\times$-seminorm at each Galerkin iteration. (Right-Bottom) The progress of the relative error in the broken $\nabla\times$-seminorm within each Galerkin iteration. }
\label{3dvariamaxcompare}
\end{figure}

Table \ref{3dvariamaxtable} shows computational complexity of the DGNN methods with one-hidden layer and two-hidden layers, respectively, for different $h$ and neurons. Figure \ref{3dvariamaxcompare} shows the relative errors of ${\bf E}-{\bf E}_r$ in the $L^2$-norm and the broken $\nabla\times$-seminorm at each Galerkin iteration, respectively. We also provide the analogous results after each training epoch. %Set the number of training epoches as 2 in {\bf Algorithm 4.1}.

\begin{center}
       \tabcaption{}%\vskip -0.3in
\label{3dvariamaxtable}
      Comparisons of computational complexity.  \vskip 0.1in
\begin{tabular}{|c|c|c|c|} \hline
   \text{Method} & \text{Neurons} & \text{DOFs} & \text{Relative} $L^2$ \text{Error} \\ \hline
  \text{Two-layers case}, $\omega h=\frac{\pi}{3}$ & 228096 &  145152 & 3.15e-4  \\ \hline
  %\text{One-layer}, $\omega h=\frac{\pi}{3}$  & 5616 & $\textcolor{red}{3\times}$ 5616 & 1.77e-3 \\  \hline
  \text{One-layer case}, $\omega h=\frac{\pi}{4}$ &  254016 &  254016 &  7.78e-4 \\  \hline
   \end{tabular}
     \end{center}

Numerical results validate that, the resulting approximate solutions generated by the DGNN method employing two-hidden-layer neural network architectures can reach the higher accuracy. Moreover, the cost of computing linear coefficients ${\bf W}^{(k)}_3$ by Galerkin least squares (\ref{dresivari}) and updating nonlinear weights ${\bf W}^{(k)}_1, {\bf W}^{(k)}_2$ by the  gradient descent algorithm (\ref{gradientIMPLE}) is smaller.

\subsection{Time-dependent linear wave equations}
Consider the second order time-dependent scalar wave equation with the source term $f$ and the coefficient $c$: %variable (see \cite{yuan2})
 \begin{equation} \label{scalarmodel}
\left\{ \begin{aligned}
       & -\Delta u  + c^{-2} \frac{\partial^2 u}{\partial t^2}=f & \text{in} \quad Q, \\
      &  u(\cdot,0)= U_0, \quad \frac{\partial u}{\partial t}(\cdot,0) = v_0 & \text{on} \quad \Omega, \\
            &     \frac{\partial u}{\partial t}  = g_D & \text{on} \quad \Gamma_D\times[0,T],
       \\
      &   -\nabla u \cdot {\bf n}^{\bf x}_{\Omega} = g_N & \text{on} \quad \Gamma_N\times[0,T].
                          \end{aligned} \right.
                          \end{equation}
Here $c$ is the wavespeed. We would like to point out that this model is different from that studied in \cite{imbertmosto} and \cite{yuan2}, attributed to the presence of source terms and variable coefficients. The purpose is to test the effectiveness of the proposed method on actual models.

Assume the space-time grid \( {\cal T}_h\) is a Cartesian-product mesh: on an internal face $F=\partial K_1\bigcap \partial K_2$, either
\begin{equation} \label{spacetimemesh}
\left\{ \begin{aligned}
     &   {\bf n}_F^{\bf x}=0  & \text{and} ~F ~\text{is called ``space-like" face, or} \\
      &  {\bf n}_F^t=0  & \text{and} ~F ~\text{is called ``time-like" face}, \\
                          \end{aligned} \right.
                          \end{equation}
where $({\bf n}_F^{\bf x},{\bf n}_F^t)$ is a unit vector of the face $F$. We denote by $\mathcal{F}_h^{\text{space}}$ the union of the internal space-like faces,  and by $\mathcal{F}_h^{\text{time}}$ the union of the internal time-like faces, respectively. Set $\mathcal{F}_h^D = \mathcal{F}_h\bigcap (\Gamma_D\times {\cal T}_{h_t}^t)$ and $\mathcal{F}_h^N = \mathcal{F}_h\bigcap (\Gamma_N\times {\cal T}_{h_t}^t)$.

The operators $\mathcal{A}, \mathcal{B}_1, \mathcal{B}_2, \mathcal{I}_1, \mathcal{I}_2, \mathcal{C}_i(i=1,\cdots,6)$ are defined by, respectively,
\beq \label{tdwave}
& \mathcal{A}= -\Delta  + c^{-2} \frac{\partial^2 }{\partial t^2}  ~\text{in}~ Q, \quad
\mathcal{B}_1 = \frac{\partial }{\partial t} ~\text{on} ~ \mathcal{F}_h^D, \quad
\mathcal{B}_2 = -\frac{\partial }{\partial {\bf n}^x_{\Omega}} ~\text{on} ~ \mathcal{F}_h^N,
\cr
& \mathcal{I}_1 = I(\cdot,0), \quad \mathcal{I}_2 = \frac{\partial }{\partial t}(\cdot,0) ~\text{on} ~ \Omega,
\cr
& \mathcal{C}_1 (u)  = \llbracket u \rrbracket_{\bf N}, ~ \mathcal{C}_2 (u)  = \llbracket \frac{\partial u}{\partial t} \rrbracket_{\bf N},
 ~ \mathcal{C}_3 (u)  = \llbracket -\nabla u \rrbracket_{\bf N}  ~\text{on}~ \mathcal{F}_h^{\text{time}},
 \cr &
 \mathcal{C}_4 (u) = \llbracket u \rrbracket_t,  ~\mathcal{C}_5 (u) = \llbracket \frac{\partial u}{\partial t} \rrbracket_t,
 ~\mathcal{C}_6 (u) = \llbracket -\nabla u \rrbracket_t ~\text{on}~ \mathcal{F}_h^{\text{space}}.
 \eq

\subsubsection{The case of constant coefficients} The wave speed is fixed at $c=10$. We test our method with the exact solution %(see \cite{yhijcom})
$$u({\bm x},t) = \text{sin}(\pi x) \text{sin}(\sqrt{2}\pi t),$$
space-time domain $Q=(0,1)\times(0,1)$, and Dirichlet boundary $\Gamma_D = \partial\Omega$.

The relaxation parameters $\{\lambda_i\}$ are initialized as follows.
\be
\lambda_0= \frac{1}{\pi^2}, \lambda_2=\lambda_4=\lambda_7=\pi^2,
\en
the other relaxation parameters are set to be 1.

We employ the proposed DGTNN method with a single hidden layer introduced in section \ref{constmodel} to solve (\ref{scalarmodel}). As introduced in the previous section \ref{local_spec} for computing the initial value $u_0$, the variational problem of the nonhomogeneous local governing equation \((\ref{cha4.n1})\)  is to find \(
u^{(1)}_k\in H^1_{IB}(Q^{\ast}_k)\) such that
\begin{equation} \label{TDconstacous} \left\{
\begin{aligned}
     &\int_{Q^{\ast}_k}(\nabla u^{(1)}_k\cdot\nabla {v}_k  - c^{-2} \frac{\partial{u^{(1)}_k}}{\partial t} \frac{\partial{v_k}}{\partial t} ) dV
     + \int_{\partial\Omega^{\ast}_q\times I_p} -\nabla u^{(1)}_k \cdot {\bf n}^{\bf x} {v}_k dS + \int_{\Omega^{\ast}_q} c^{-2} \frac{\partial{u^{(1)}_k}}{\partial t}(\cdot,t_p) v_k dS \\
     & \quad \quad\quad   =
\int_{Q^{\ast}_k}f {v}_k dV, \quad \quad\quad \forall v_k\in
H^1_{IB}(Q^{\ast}_k)~(k = 1,2,\ldots,N),
                          \end{aligned} \right.
                          \end{equation}
where \(H^1_{IB}(Q^{\ast}_k) = \{ v: v\in H^1(Q^{\ast}_k), ~v(\cdot, t_{p-1})=0 ~\text{on} ~\Omega^{\ast}_q ~\text{and}~ \frac{\partial{v}}{\partial t}=0 ~\text{on} ~ \partial\Omega^{\ast}_q\times I_p  \}.\)

Local spectral elements of order $m=5$  are chosen to represent the initial value $u_0$ in the DGTNN iterative algorithm.
 The discrete variational problems of (\ref{TDconstacous}) are: to find $u^{(1)}_{k,h}\in S_5(Q^{\ast}_k) \cap H^1_{IB}(Q^{\ast}_k)$ such that
\begin{equation} \label{TDconstacousvari} \left\{
\begin{aligned}
     &\int_{Q^{\ast}_k}( \nabla u^{(1)}_{k,h}\cdot\nabla {v}_k  - c^{-2} \frac{\partial{u^{(1)}_{k,h}}}{\partial t} \frac{\partial{v_k}}{\partial t} ) dV
     + \int_{\partial\Omega^{\ast}_q\times I_p} -\nabla u^{(1)}_{k,h} \cdot {\bf n}^{\bf x} {v}_k dS + \int_{\Omega^{\ast}_q} c^{-2} \frac{\partial{u^{(1)}_{k,h}}}{\partial t}(\cdot,t_p) v_k dS \\
     & \quad \quad\quad   =
\int_{Q^{\ast}_k}f {v}_k dV,\quad \quad\quad \forall v_k\in
S_5(Q^{\ast}_k)  \cap H^1_{IB}(Q^{\ast}_k) ~(k = 1,2,\ldots,N).
                          \end{aligned} \right.
                          \end{equation}

Subsequently, we need to construct the Trefftz neural networks.
First, one can see the fact that, for any polynomial function $Q$ and any unit vector ${\bf d}\in \mathbb{R}^d$, the space-time field $Q ({\bf d}\cdot {\bm x}-ct)$ is solution of the second order wave equation (\ref{scalarmodel}) without initial and boundary conditions. Then for each $l=0, \cdots, p$, and a fixed monomial $Q_l \in \mathbb{P}^l(\mathbb{R})$ of degree exactly $l$, define a {\it discontinuous} Trefftz neural network function $\varphi^{\theta}: \mathbb{R}^{d+1} \rightarrow \mathbb{C}$ as follows:
\be \label{nn_tdacousconst_general2}
\varphi^{\theta}({\bf x})|_{\Omega_k} = \sum_{l=0}^p\sum_{j=1}^{d_l} c_{l,j}^{(k)} Q_l({\bf d}^{(k)}_{l,j}\cdot {\bm x}-ct),
\en
where $\{{\bf d}_{l,j}^{(k)}\in \mathbb{R}^d: |{\bf d}_{l,j}^{(k)}|=1\}$ are different propagation directions on
$k$-th element $\Omega_k\in {\cal T}_h$; $c_{l,j}^{(k)} \in \mathbb{C}$ are element-wise coefficients. Clearly, the width $n$ of the network is equal to $\sum_{l=0}^p d_l$. Let $\text{T}_n({\mathcal T}_h)$ be the set of all functions of the form (\ref{nn_tdacousconst_general2}). For each $l=1, \cdots, p$, the propagation directions $\{{\bf d}_{l,j}^{(k)}, j=1,\cdots,d_l\}$ are initialized to be uniformly distributed on the unit circle.

We also employ the DGNN method with two hidden layers introduced in section \ref{twolayer_sec} to solve
(\ref{scalarmodel}). %Employ two hidden layers, and sigmoid activation function
Initialize nonlinear parameters as follows. Weights ${\bf W}^{(k)}_1\in \mathbb{C}^{m_1\times 2}$ in the first hidden layer are uniformly distributed on the unit circle.
\be\nonumber
{\bf W}^{(k)}_1(j,:)=(cos\theta_j, sin\theta_j), ~~~\theta_j=\frac{2\pi}{m_1}j, j=1,\cdots,m_1.
\en
Set biases ${\bf b}^{(k)}_1={\bf 1}\in \mathbb{C}^{m_1}$. Weights ${\bf W}^{(k)}_2\in \mathbb{C}^{m_2\times m_1}$ in the second hidden layer are set as
${\bf W}^{(k)}_2 = -c Q$, where $m_2\geq m_1$, and the matrix $Q$ is composed of the first $m_1$ columns of a $m_2$-order random orthogonal matrix. Set biases ${\bf b}^{(k)}_2={\bf 1}\in \mathbb{C}^{m_2}$.

The other parameters are set as follows.
\be \nonumber
 \left\{\begin{array}{ll} h=\frac{1}{12}, n_r=2r+1, u_0 = u^{(1)}_h & \text{for DGTNN}\\
h=\frac{1}{12}, m_1= m_2=n_r=2r+3, u_0 = 0 & \text{for DGNN}.
\end{array}\right.
\en

\begin{figure}[H]
%\vspace{-2cm}thb
\begin{center}
\begin{tabular}{cc}
\epsfxsize=0.4\textwidth\epsffile{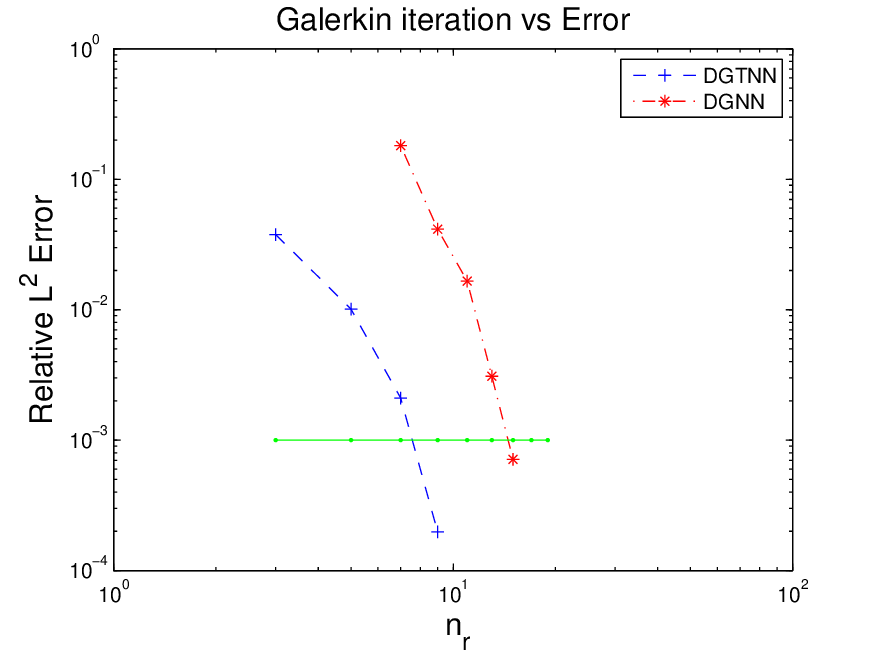}&
%COMPARE_2dhelm_galer_l2.eps
\epsfxsize=0.4\textwidth\epsffile{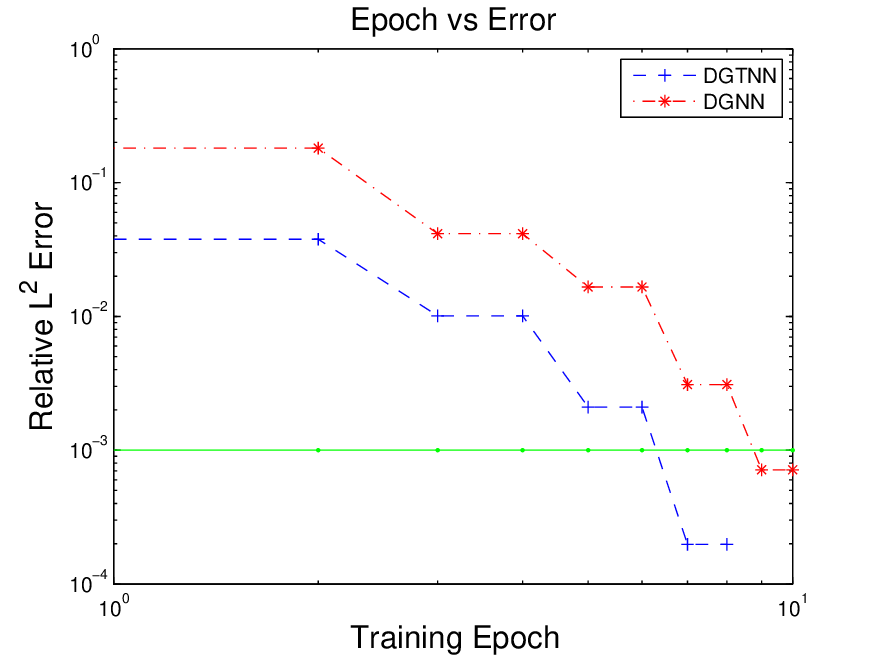}\\
\epsfxsize=0.4\textwidth\epsffile{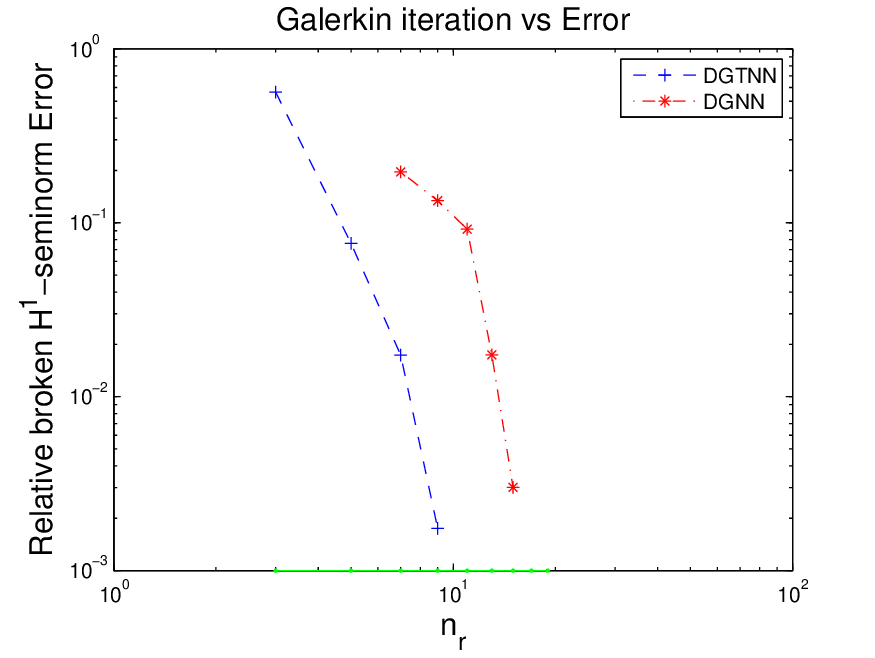}&
\epsfxsize=0.4\textwidth\epsffile{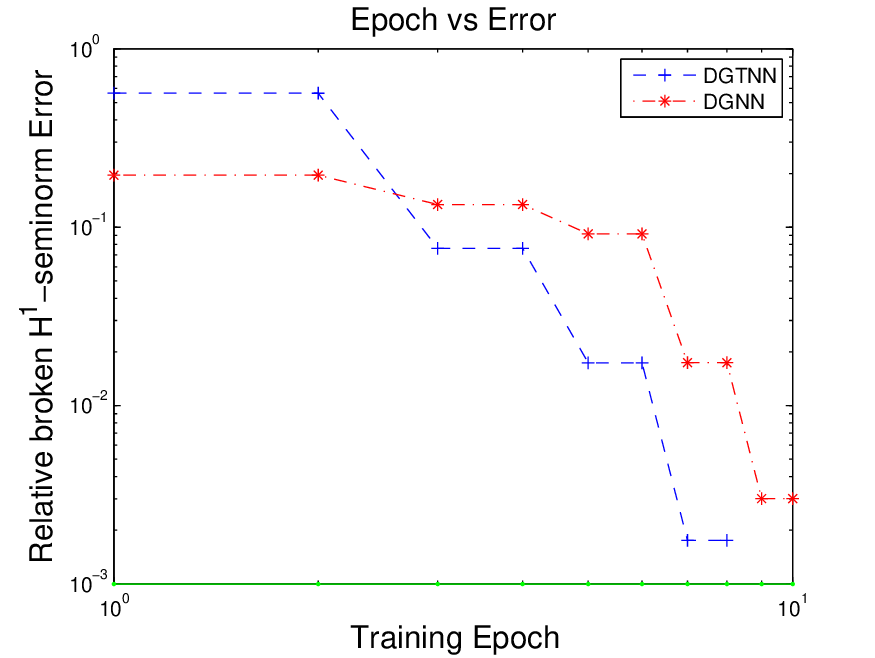}\\
\end{tabular}
\end{center}
 \caption{Displacement of a string (Time-dependent linear wave equations with constant coefficients). (Left-Up) Relative error in the $L^2$-norm at each Galerkin iteration. (Right-Up) The progress of the relative error in the $L^2$-norm within each Galerkin iteration. The x-axis thus denotes the cumulative training epoch over all Galerkin iterations. (Left-Bottom) Relative error in the broken $H^1$-seminorm at each Galerkin iteration. (Right-Bottom) The progress of the relative error in the broken $H^1$-seminorm within each Galerkin iteration. }
\label{3dconsTDaoucticwavenn1}
\end{figure}

Table \ref{3dconsTDACOUStable} shows computational complexity of the DGTNN method and DGNN method for different $h$ and neurons.
Figure \ref{3dconsTDaoucticwavenn1} shows the relative errors of $u-u_r$ in the $L^2$-norm and the broken $H^1$-seminorm at each Galerkin iteration, respectively. We also provide the analogous results after each training epoch.

\begin{center}
       \tabcaption{}%\vskip -0.3in
\label{3dconsTDACOUStable}
      Comparisons of computational complexity.  \vskip 0.1in
\begin{tabular}{|c|c|c|c|} \hline
   \text{Method} & \text{Neurons} & \text{DOFs} & \text{Relative} $L^2$ \text{Error} \\ \hline
  \text{DGTNN} & 1296 &  1296 & 1.98e-4  \\ \hline
  \text{DGNN} &  4320 &  2160  & 7.13e-4 \\  \hline
   \end{tabular}
     \end{center}

Compared with the DGNN method with general activation functions, Trefftz neural networks can achieve higher accuracy with smaller computational costs, that is, relatively coarser spatial width and smaller number of neurons. %Particularly, Trefftz neural networks only require $h$

\subsubsection{The case of variable coefficients}
Consider the variable wavespeed parameter $c=x+1$.
We still test our method with the exact solution %(see \cite{yhijcom})
$$u({\bm x},t) = \text{sin}(\pi x) \text{sin}(\sqrt{2}\pi t),$$
space-time domain $Q=(0,1)\times(0,1)$, and Dirichlet boundary $\Gamma_D = \partial\Omega$.

The DGNN method with two hidden layers introduced in the last subsection is employed. Set $u_0=0$, $h=\frac{1}{16}$, $m_1 = m_2 = n_r=2r+3$ for each Galerkin iteration $r$.

\begin{figure}[H]
%\vspace{-2cm}thb
\begin{center}
\begin{tabular}{cc}
\epsfxsize=0.4\textwidth\epsffile{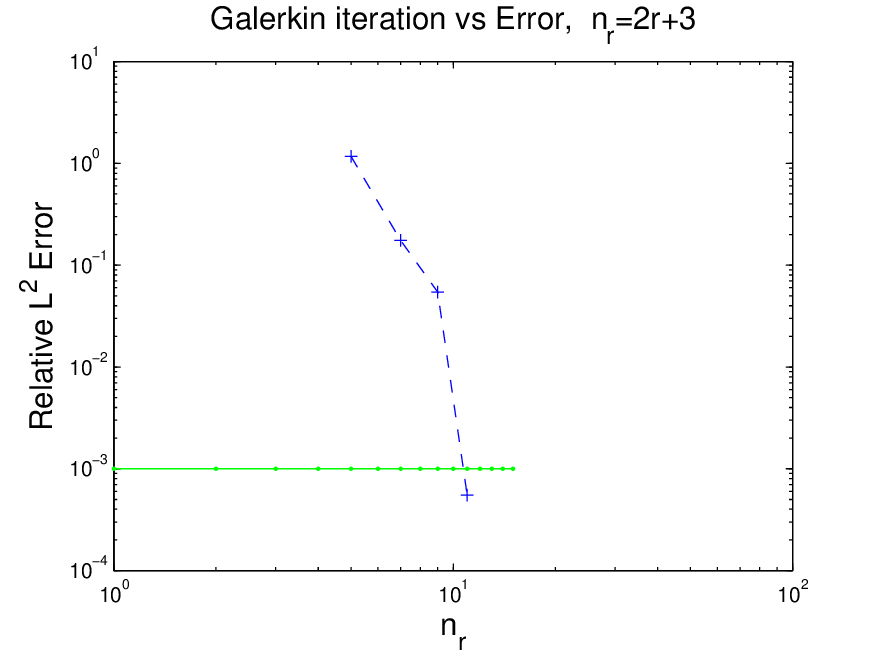}&
\epsfxsize=0.4\textwidth\epsffile{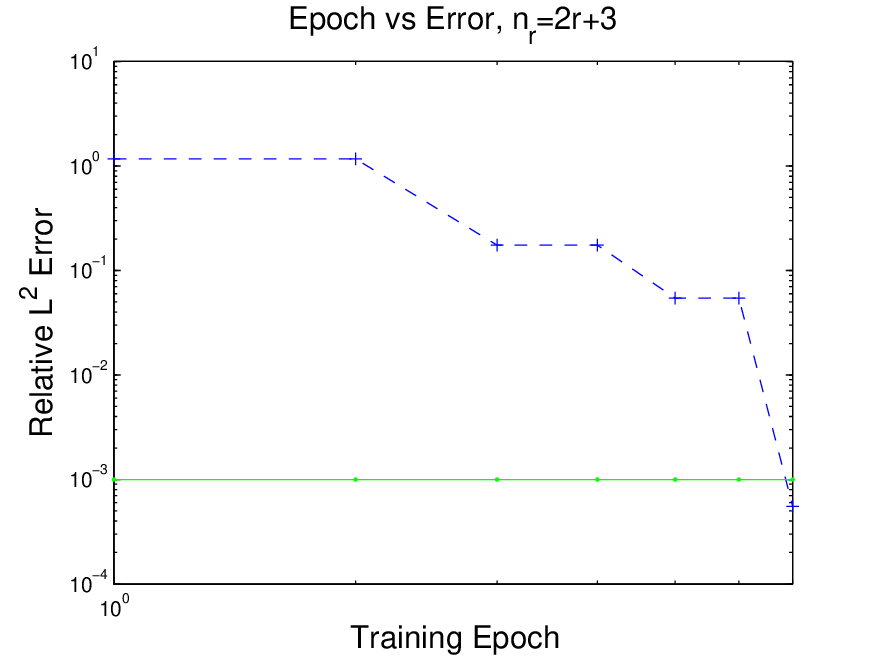}\\
\epsfxsize=0.4\textwidth\epsffile{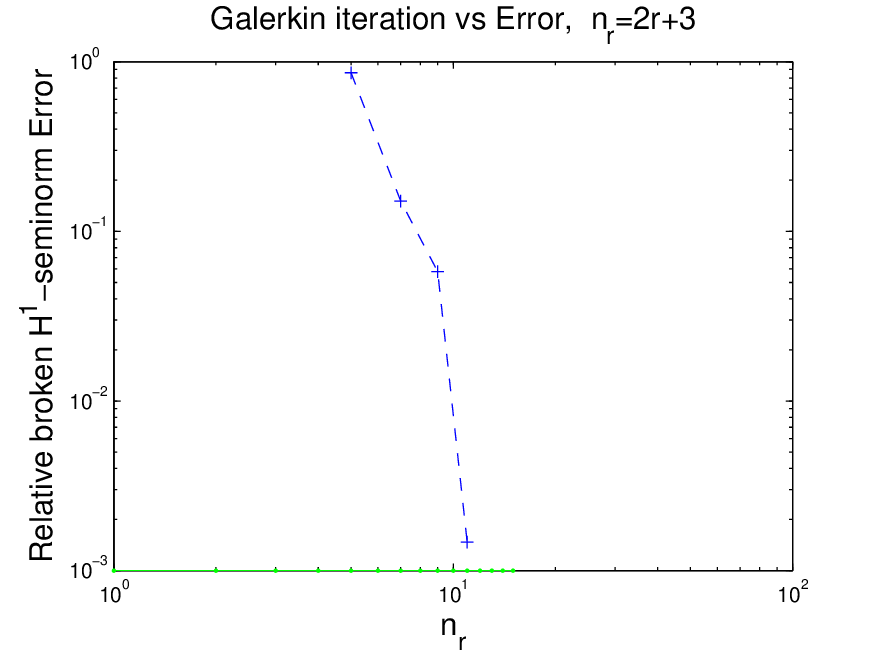}&
\epsfxsize=0.4\textwidth\epsffile{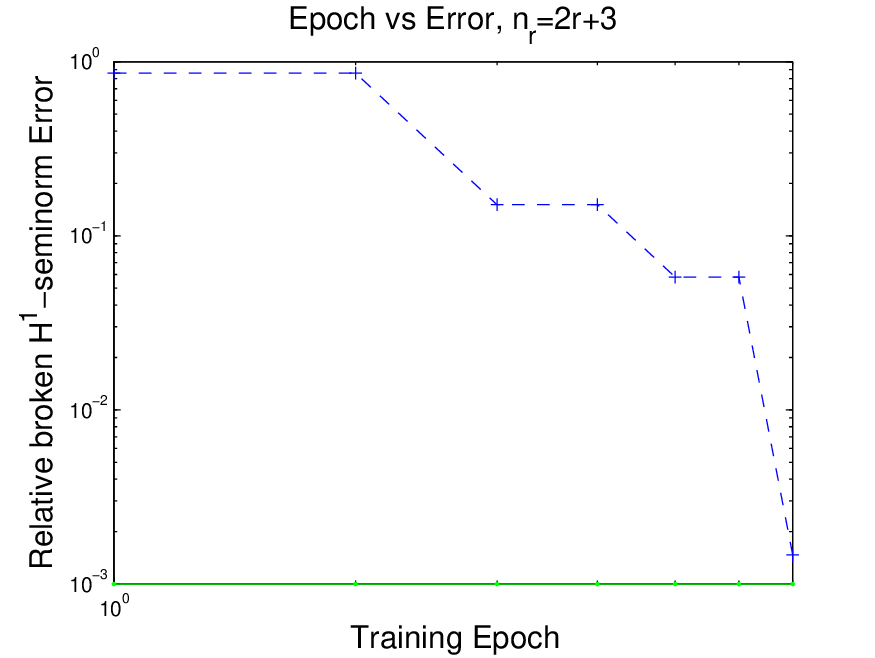}\\
\end{tabular}
\end{center}
 \caption{Displacement of a string (Time-dependent acoustic wave equation with variable coefficients). (Left-Up) Relative error in the $L^2$-norm at each Galerkin iteration. (Right-Up) The progress of the relative error in the $L^2$-norm within each Galerkin iteration. The x-axis thus denotes the cumulative training epoch over all Galerkin iterations. (Left-Bottom) Relative error in the broken $H^1$-seminorm at each Galerkin iteration. (Right-Bottom) The progress of the relative error in the broken $H^1$-seminorm within each Galerkin iteration. }
\label{2dvariTDacousticwave_nn1}
\end{figure}

Table \ref{2dvariTDtable} shows computational complexity of the DGNN method. Figure \ref{2dvariTDacousticwave_nn1} shows the relative errors of $u-u_r$ in the $L^2$-norm and the broken $H^1$-seminorm at each Galerkin iteration, respectively. We also provide the analogous results after each training epoch.

\begin{center}
       \tabcaption{}%\vskip -0.3in
\label{2dvariTDtable}
      Computational complexity of the DGNN method.  \vskip 0.1in
\begin{tabular}{|c|c|c|} \hline
    \text{Neurons} & \text{DOFs} & \text{Relative} $L^2$ \text{Error} \\ \hline
   5632 & 2816 &  5.51e-4 \\ \hline
   \end{tabular}
     \end{center}

Numerical results validate that, the resulting approximate solutions generated by the DGNN method employing two-hidden-layer neural network architectures can reach the given accuracy.

\subsection{Time-dependent electromagnetic wave equations}
 We consider the second order electromagnetic wave IBVP posed on
 a space-time domain $Q=\Omega\times I$, where $\Omega \subset \mathbb{R}^3$ is an open bounded Lipschitz polytope and $I=(0,T),T>0$.
 ${\bf n}_{\Omega}$ is an outward-pointing unit normal vector on $\Gamma=\partial\Omega$. The second order model reads as (see \cite{ALZOU,BZZOU})
  \begin{equation} \label{secondmodel}
\left\{ \begin{aligned}
     &   \varepsilon \frac{\partial^2 {\bf E}}{\partial t^2} + \nabla \times ( \mu^{-1}\nabla \times {\bf E} ) ={\bf f} & \text{in} \quad Q, \\
      &  {\bf E}(\cdot,0)={\bf F}_0, ~~\frac{\partial{\bf E}}{{\partial t}}(\cdot,0)={\bf H}_0 & \text{on} \quad \Omega,
        \\
      &  {\bf n}_{\Omega}\times  {\bf E} = {\bf g} & \text{on} \quad \Gamma \times[0,T].
                          \end{aligned} \right.
                          \end{equation}
Here  ${\bf E}=(E_x,E_y,E_z)^T$, ${\bf H}=(H_x,H_y,H_z)^T$; ${\bf F}_0 \in {\bf H}_0(\text{curl};\Omega)$ and ${\bf H}_0 \in {\bf L}^2(\Omega)$ are the given source data; \( {\bf g}\in {\bf L}^2_{\bf T}(\Gamma\times[0,T])\). The permittivity \(\varepsilon\) and the permeability \(\mu\) are piecewise constant satisfying $0<c_{\varepsilon}\leq \varepsilon \leq C_{\varepsilon}$ and $0<c_{\mu}\leq \mu \leq C_{\mu}$.

The operators $\mathcal{A}, \mathcal{B}_1, \mathcal{I}_1, \mathcal{I}_2, \mathcal{C}_i(i=1,\cdots,6)$ are defined by, respectively,
\beq \label{tdmax}
& \mathcal{A}= \varepsilon \frac{\partial^2 }{\partial t^2} +  \nabla \times ( \mu^{-1}\nabla \times \cdot )  ~\text{in}~ Q, \quad
\mathcal{B}_1 = {\bf n}_{\Omega}\times \cdot ~\text{on} ~ \mathcal{F}_h^D,
\cr
& \mathcal{I}_1 = I(\cdot,0), \quad \mathcal{I}_2 = \frac{\partial }{\partial t}(\cdot,0) ~\text{on} ~ \Omega,
\cr
& \mathcal{C}_1 ({\bf E})  = \llbracket {\bf E} \rrbracket_{\bf T}, ~ \mathcal{C}_2 ({\bf E})  = \llbracket \frac{\partial{\bf E}}{{\partial t}} \rrbracket_{\bf T},
 ~ \mathcal{C}_3 ({\bf E})  = \llbracket \nabla \times {\bf E} \rrbracket_{\bf T}  ~\text{on}~ \mathcal{F}_h^{\text{time}},
 \cr &
 \mathcal{C}_4 ({\bf E}) = \llbracket {\bf E} \rrbracket_t,  ~\mathcal{C}_5 ({\bf E}) = \llbracket \frac{\partial{\bf E}}{{\partial t} }\rrbracket_t,
 ~\mathcal{C}_6 ({\bf E}) = \llbracket \nabla \times {\bf E} \rrbracket_t ~\text{on}~ \mathcal{F}_h^{\text{space}}.
 \eq

The relaxation parameters $\{\lambda_i\}$ are initialized as follows.
\be
\lambda_0= \frac{1}{\pi^2}, \lambda_1=\lambda_2=\lambda_4=\lambda_7=\pi^2,
\en
the other relaxation parameters are set to be 1.

\subsubsection{The case of constant coefficients}

We test our method with the material coefficients $\mu =1, \varepsilon = 1$. Consider the following analytical solution
\be \label{TDmaxequation}
{\bf E}_{\text{ex}} =
\left(
  \begin{array}{c}
    sin(\pi  y) sin(\pi  z) cos(\pi t) \\
    sin(\pi  x) sin(\pi  z) sin(\pi t) \\
    cos(\pi  x) sin(\pi  y) cos(\pi t) \\
  \end{array}
\right)
\en
in a space-time domain $Q=(0,1)^3\times(0,1)$.

We employ the proposed DGTNN method with a single hidden layer introduced in section \ref{constmodel} to solve (\ref{secondmodel}). As introduced in the previous section \ref{local_spec} for computing the initial value ${\bf E}_0$, the variational problem of the nonhomogeneous local governing equation \((\ref{cha4.n1})\)  is to find \(
{\bf E}^{(1)}_k\in {\bf H}_{IB}(Q^{\ast}_k)\) such that
\begin{equation} \label{TDconstmax} \left\{
\begin{aligned}
     &\int_{Q^{\ast}_k}( -\varepsilon \frac{\partial{\bf E}^{(1)}_k}{{\partial t}} \cdot \frac{\partial{\bf F}_k}{{\partial t}} + \mu^{-1}\nabla \times {\bf E}^{(1)}_k  \cdot \nabla \times {\bf F}_k ) dV
     + \int_{\Omega^{\ast}_q} \varepsilon \frac{\partial{\bf E}^{(1)}_k(\cdot,t_p)}{{\partial t}} \cdot  {\bf F}_k(\cdot,t_p) dS \\
     & \quad \quad\quad   =
\int_{Q^{\ast}_k}{\bf f}\cdot {\bf F}_k dV, \quad \quad\quad \forall {\bf F}_k\in
{\bf H}_{IB}(Q^{\ast}_k)~(k = 1,2,\ldots,N),
                          \end{aligned} \right.
                          \end{equation}
where \({\bf H}_{IB}(Q^{\ast}_k) = \{ {\bf F}: {\bf F} \in { L}^2(I_p; {\bf H}(\text{curl};\Omega^{\ast}_q)), ~ \dot{{\bf F}}
\in { L}^2(I_p; {\bf L}^2(\Omega^{\ast}_q)), ~{\bf F}(\cdot,t_{p-1})={\bf 0}  \}.\)

Local spectral elements of order $m=5$  are chosen to represent the initial value ${\bf E}_0$ in the DGTNN iterative algorithm.
The discrete variational problems of (\ref{TDconstmax}) are: to find ${\bf E}^{(1)}_{k,h}\in {\bf S}_5(Q^{\ast}_k) \cap {\bf H}_{IB}(Q^{\ast}_k)$ such that
\begin{equation} \label{TDconstmaxvari} \left\{
\begin{aligned}
     &\int_{Q^{\ast}_k}( -\varepsilon \frac{\partial{\bf E}^{(1)}_{k,h}}{{\partial t}} \cdot \frac{\partial{\bf F}_k}{{\partial t}} + \mu^{-1}\nabla \times {\bf E}^{(1)}_{k,h}  \cdot \nabla \times {\bf F}_k ) dV
     + \int_{\Omega^{\ast}_q} \varepsilon \frac{\partial{\bf E}^{(1)}_{k,h}(\cdot,t_p)}{{\partial t}} \cdot  {\bf F}_k(\cdot,t_p) dS \\
     & \quad \quad\quad   =
\int_{Q^{\ast}_k}{\bf f}\cdot {\bf F}_k dV, \quad \quad\quad \forall {\bf F}_k\in
{\bf S}_5(Q^{\ast}_k)  \cap {\bf H}_{IB}(Q^{\ast}_k)~(k = 1,2,\ldots,N).
                          \end{aligned} \right.
                          \end{equation}

Later, we need to construct the Trefftz neural networks.
First, one can see the fact that, for any polynomial function $Q$ and any unit propagation vector ${\bf d}\in \mathbb{R}^d$, the space-time vector field ${\bf e} Q ({\bf d}\cdot {\bm x} - c t)$ is solution of the second order wave equation (\ref{secondmodel}) without initial and boundary conditions, where a unit real polarization vector ${\bf e}$ is orthogonal to ${\bf d}$, and $c=\frac{1}{\sqrt{\varepsilon\mu}}$. Then for each $l=0, \cdots, p$, and a fixed monomial $Q_l \in \mathbb{P}^l(\mathbb{R})$ of degree exactly $l$, define a {\it discontinuous} Trefftz neural network function $\varphi^{\theta}: \mathbb{R}^{d+1} \rightarrow \mathbb{C}$ as follows:
\be \label{nn_tdmaxconst_general2}
\varphi^{\theta}({\bf x})|_{\Omega_k} = \sum_{l=0}^p\sum_{j=1}^{d_l} c_{l,j}^{(k)} {\bf e}_{l,j}^{(k)} Q_l({\bf d}^{(k)}_{l,j}\cdot {\bm x}-ct)
+ c_{l,j+d_l}^{(k)}  {\bf e}_{l,j+d_l}^{(k)} Q_l({\bf d}^{(k)}_{l,j}\cdot {\bm x}-ct),
\en
where $d_0=1, d_l=(l+1)^2$ for $l=1,\cdots, p$; $\{{\bf d}_{l,j}^{(k)}\in \mathbb{R}^d: |{\bf d}_{l,j}^{(k)}|=1\}$ are different propagation directions on
$k$-th element $\Omega_k\in {\cal T}_h$; $c_{l,j}^{(k)} \in \mathbb{C}$ are element-wise coefficients. For $j=1,\cdots,d_l$, every unit real polarization vector ${\bf e}_{l,j}^{(k)}$ is orthogonal to ${\bf d}_{l,j}^{(k)}$, and ${\bf e}_{l,j+d_l}^{(k)}={\bf e}_{l,j}^{(k)}\times {\bf d}_{l,j}^{(k)}$.
 Clearly, the width $n$ of the network is equal to $\sum_{l=0}^p 2d_l$. Let ${\bf T}_{n}({\mathcal T}_h)$ be the set of all functions of the form (\ref{nn_tdmaxconst_general2}). For each $l=1, \cdots, p$, the propagation directions $\{{\bf d}_{l,j}^{(k)}, j=1,\cdots,d_l\}$ are initialized to be uniformly distributed on the unit circle.

We also employ the DGNN method with two hidden layers introduced in section \ref{twolayer_sec} to solve
(\ref{secondmodel}). %the vector composed of sigmoid activation functions%Employ two hidden layers, and sigmoid activation function
Initialize nonlinear parameters as follows. Weights ${\bf W}^{(k)}_1 \in \mathbb{C}^{m_1\times 4}$ in the first hidden layer are decomposed into two parts ${\bf W}^{(k)}_1 = ({\bf W}^{(k)}_{1,s}, {\bf W}^{(k)}_{1,t})$, where ${\bf W}^{(k)}_{1,s}\in \mathbb{C}^{m_1\times 3}$ is generated by the spherical codes from \cite{refsite} and ${\bf W}^{(k)}_{1,t}=-c \cdot{\bf 1}\in \mathbb{C}^{m_1\times 1}$. Set biases ${\bf b}^{(k)}_1={\bf 1}\in \mathbb{C}^{m_1}$. Weights ${\bf W}^{(k)}_2\in \mathbb{C}^{m_2\times m_1}$ in the second hidden layer are set as
${\bf W}^{(k)}_2 = Q$, where $m_2\geq m_1$, the matrix $Q$ is composed of the first $m_1$ columns of a $m_2$-order random orthogonal matrix. Set biases ${\bf b}^{(k)}_2={\bf 1}\in \mathbb{C}^{m_2}$. %The reason why the nonlinear parameters ${\bf W}^{(k)}_i$ ($i=1,2$) are initialized in this model-driven manner can refer to section \ref{variahelm}.

The other parameters are set as follows.
\be \nonumber
 \left\{\begin{array}{ll} h=\frac{1}{12}, n_r=\sum\limits_{l=0}^r 2d_l, {\bf E}_0= {\bf E}^{(1)}_h  & \text{for DGTNN}\\
h=\frac{1}{12}, m_1=m_2=n_r=(r+3)^2, {\bf E}_0= {\bf 0}  & \text{for DGNN}.
\end{array}\right.
\en

\begin{figure}[H]
%\vspace{-2cm}thb
\begin{center}
\begin{tabular}{cc}
\epsfxsize=0.4\textwidth\epsffile{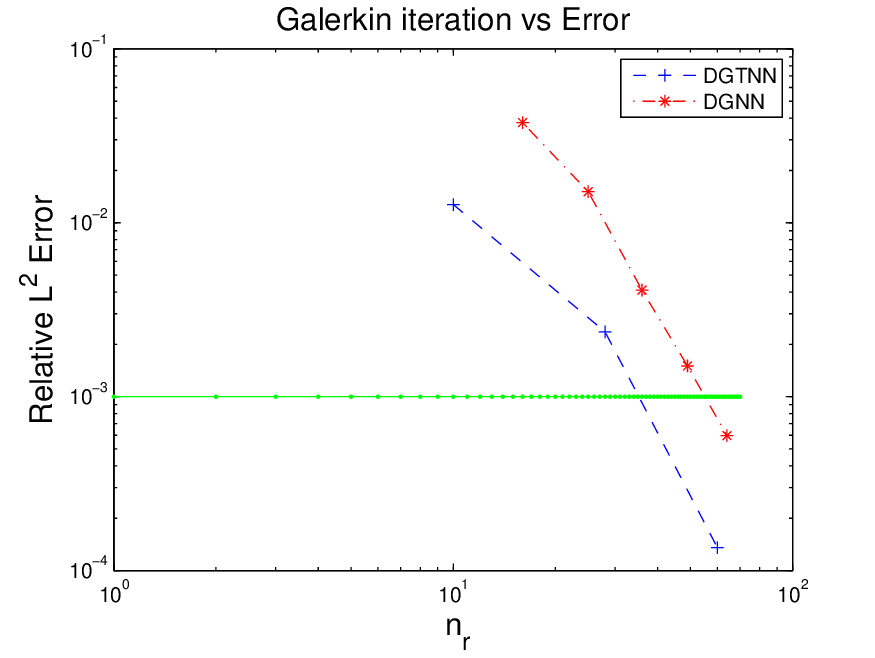}&
%COMPARE_2dhelm_galer_l2.eps
\epsfxsize=0.4\textwidth\epsffile{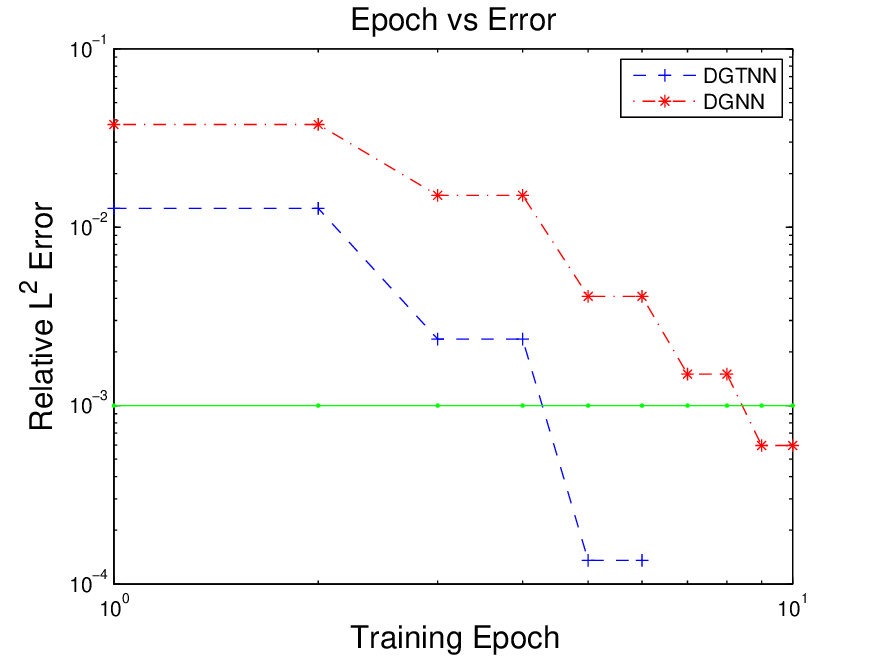}\\
\epsfxsize=0.4\textwidth\epsffile{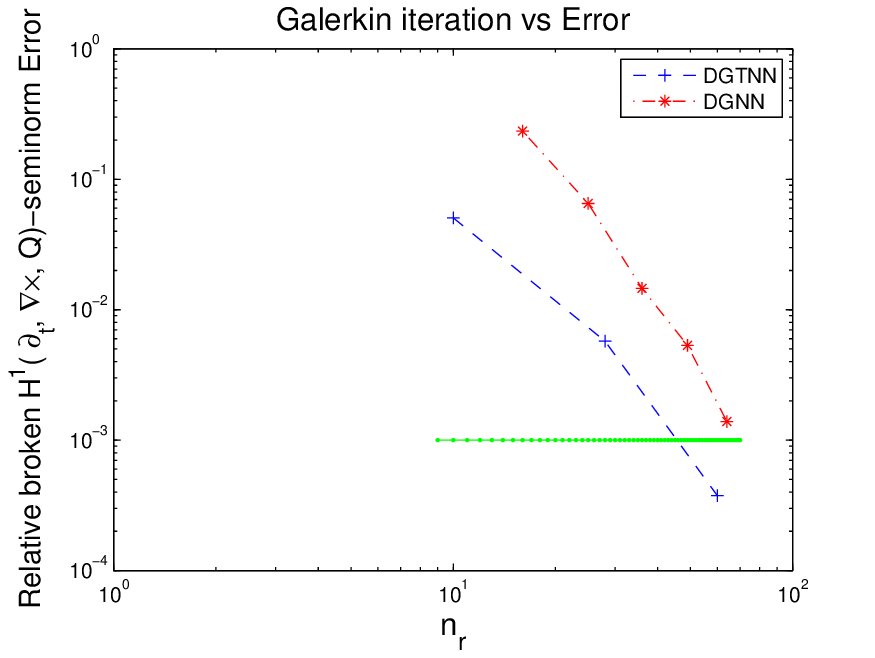}&
\epsfxsize=0.4\textwidth\epsffile{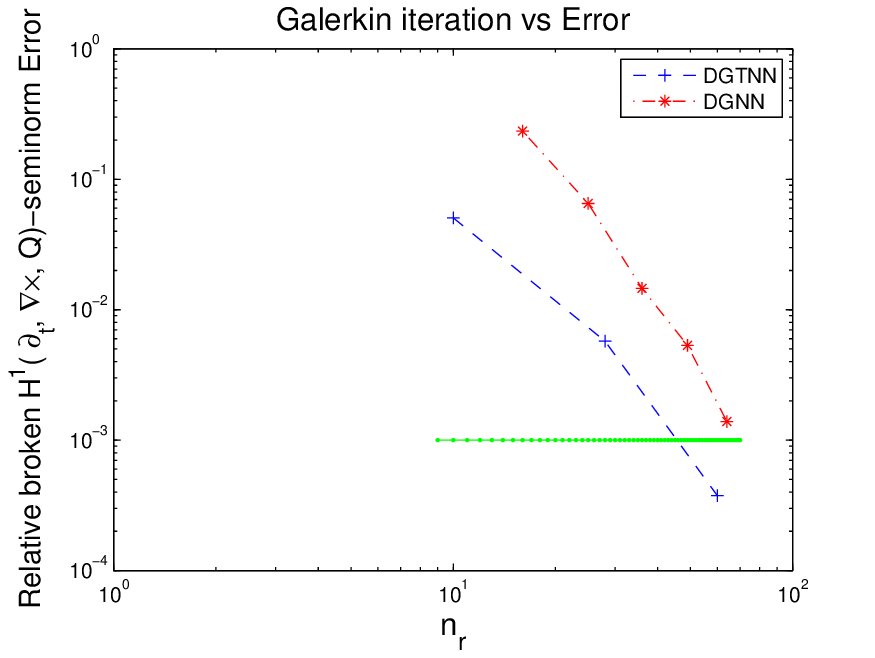}\\
\end{tabular}
\end{center}
 \caption{Displacement of a string (Maxwell's equations in three dimension). (Left-Up) Relative error in the $L^2$-norm at each Galerkin iteration. (Right-Up) The progress of the relative error in the $L^2$-norm within each Galerkin iteration. The x-axis thus denotes the cumulative training epoch over all Galerkin iterations. (Left-Bottom) Relative error in the broken $H^1(\frac{\partial}{\partial t}, \nabla\times, Q)$-seminorm at each Galerkin iteration. (Right-Bottom) The progress of the relative error in the broken $H^1(\frac{\partial}{\partial t}, \nabla\times, Q)$-seminorm within each Galerkin iteration. }
\label{3dconsTDmaxnn1}
\end{figure}

Table \ref{3dconsTDmaxtable} shows computational complexity of the DGTNN method and DGNN method for different $h$ and neurons.
Figure \ref{3dconsTDmaxnn1} shows the relative errors of ${\bf E}-{\bf E}_r$ in the $L^2$-norm and the broken $H^1(\frac{\partial}{\partial t}, \nabla\times, Q)$-seminorm at each Galerkin iteration, respectively, where $|{\bf F}|^2_{H^1(\frac{\partial}{\partial t}, \nabla\times, Q)}=(||\frac{\partial{\bf F}}{\partial t}||^2_{L^2(Q)}+||\nabla\times{\bf F}||^2_{L^2(Q)})^{\frac{1}{2}}$. We also provide the analogous results after each training epoch. %Set the number of training epoches as 2 in {\bf Algorithm 4.1}.

\begin{center}
       \tabcaption{}%\vskip -0.3in
\label{3dconsTDmaxtable}
      Comparisons of computational complexity.  \vskip 0.1in
\begin{tabular}{|c|c|c|c|} \hline
   \text{Method} & \text{Neurons} & \text{DOFs} & \text{Relative} $L^2$ \text{Error} \\ \hline
  \text{DGTNN} & 103680 &  103680 & 1.36e-4  \\ \hline
  \text{DGNN} &  221184 &   110592 &  5.98e-4 \\  \hline
   \end{tabular}
     \end{center}

Compared with the DGNN method with general activation functions, Trefftz neural networks can achieve higher accuracy with smaller computational costs, that is, relatively coarser spatial width and smaller number of neurons. %Particularly, Trefftz neural networks only require $h$

\subsubsection{The case of variable coefficients}
We still consider the analytical solution (\ref{TDmaxequation}), but with the different material coefficient $\varepsilon = 1 + x+y+z$.

The DGNN method with two hidden layers introduced in the last subsection is employed.  %Employ two hidden layers, and sigmoid activation function
 Set $h=\frac{1}{16}$, $m_1 = m_2 = n_r=(r+3)^2$ for each Galerkin iteration $r$.

\begin{figure}[H]
%\vspace{-2cm}thb
\begin{center}
\begin{tabular}{cc}
\epsfxsize=0.4\textwidth\epsffile{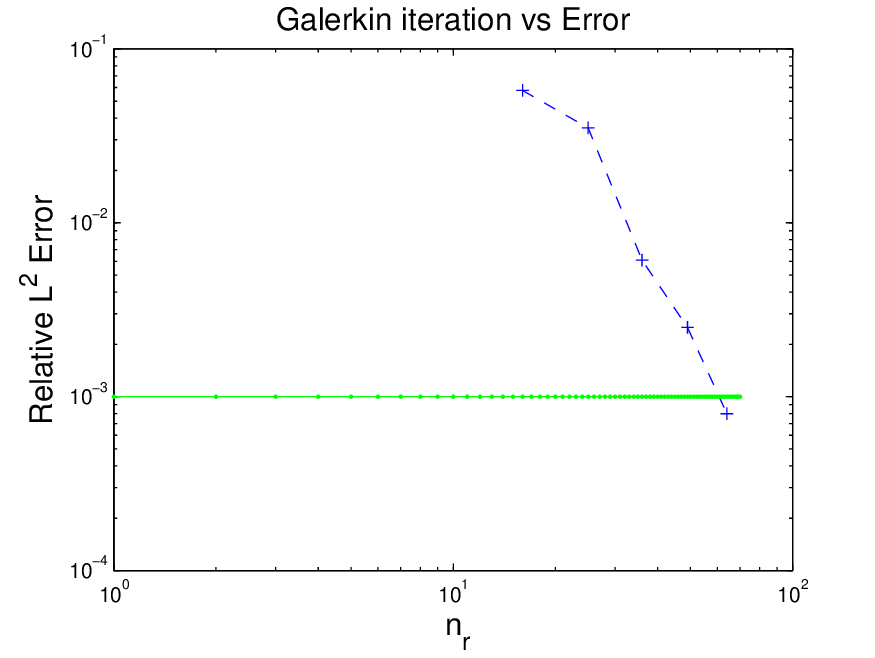}&
\epsfxsize=0.4\textwidth\epsffile{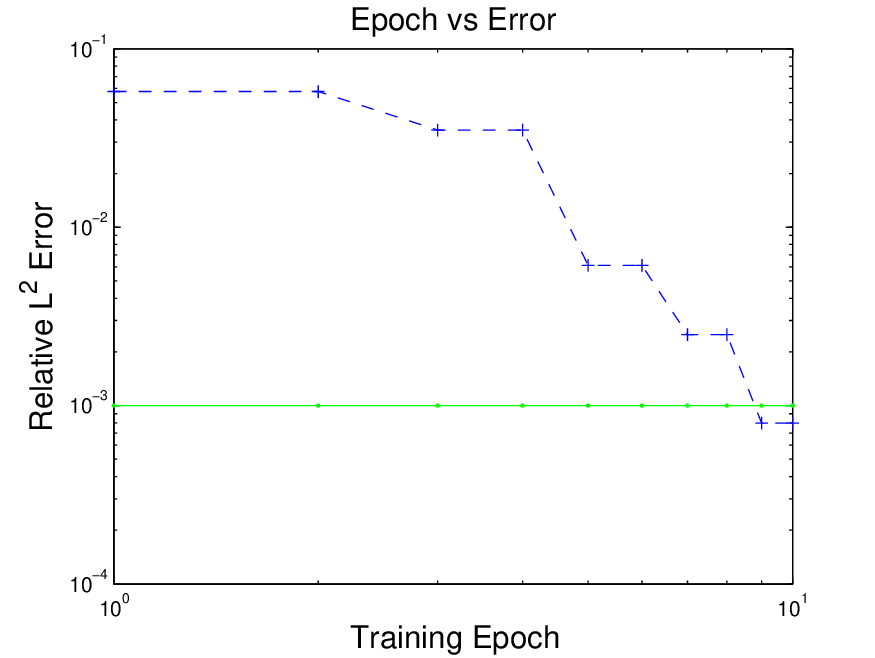}\\
\epsfxsize=0.4\textwidth\epsffile{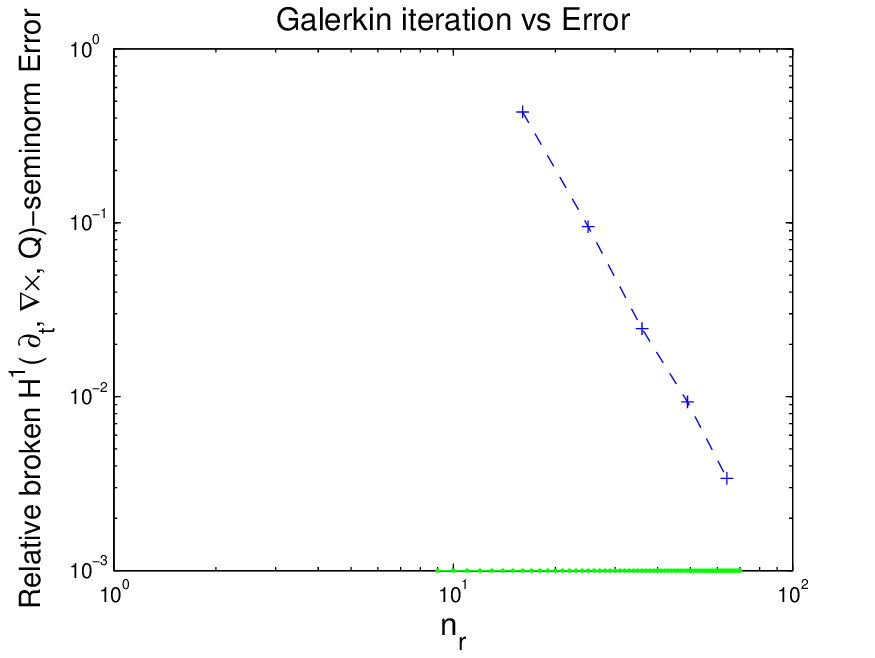}&
\epsfxsize=0.4\textwidth\epsffile{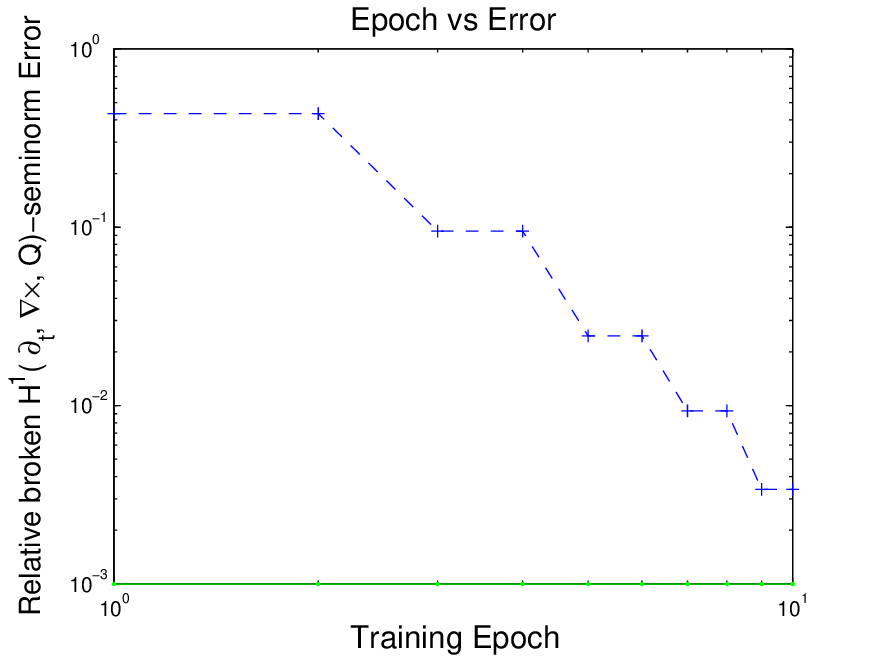}\\
\end{tabular}
\end{center}
 \caption{Displacement of a string (Maxwell's equations in three dimension). (Left-Up) Relative error in the $L^2$-norm at each Galerkin iteration. (Right-Up) The progress of the relative error in the $L^2$-norm within each Galerkin iteration. The x-axis thus denotes the cumulative training epoch over all Galerkin iterations. (Left-Bottom) Relative error in the broken $H^1(\frac{\partial}{\partial t}, \nabla\times, Q)$-seminorm at each Galerkin iteration. (Right-Bottom) The progress of the relative error in the broken $H^1(\frac{\partial}{\partial t}, \nabla\times, Q)$-seminorm within each Galerkin iteration. }
\label{2dvariTDmax_nn1}
\end{figure}

Table \ref{2dvariTDmaxtable} shows computational complexity of the DGNN method. Figure \ref{2dvariTDmax_nn1} shows the relative errors of ${\bf E}-{\bf E}_r$ in the $L^2$-norm and the broken $H^1(\frac{\partial}{\partial t}, \nabla\times, Q)$-seminorm at each Galerkin iteration, respectively. We also provide the analogous results after each training epoch. %Set the number of training epoches as 2 in {\bf Algorithm 4.1}.

\begin{center}
       \tabcaption{}%\vskip -0.3in
\label{2dvariTDmaxtable}
      Computational complexity of the DGNN method.  \vskip 0.1in
\begin{tabular}{|c|c|c|} \hline
    \text{Neurons} & \text{DOFs} & \text{Relative} $L^2$ \text{Error} \\ \hline
   524288 & 262144 &   7.89e-4 \\ \hline
   \end{tabular}
     \end{center}

Numerical results validate that, the resulting approximate solutions generated by the DGNN method employing two-hidden-layer neural network architectures can reach the given accuracy.

\section{Conclusion}
In this paper we have introduced a framework of the DGNN method for iteratively approximating residuals. Within this framework, the desired neural network approximate solution is recursively supplemented by solving a sequence of quasi-minimization problems associated with the underlying loss functionals and the adaptively augmented discontinuous neural network sets without the assumption on the boundedness of the neural network parameters. We further propose a DGTNN method with a single hidden layer to reduce the computational costs and storage requirements. Numerical experiments confirm that compared to existing PINN algorithms, the proposed DGNN method is able to improve the relative $L^2$ error by at least one order of magnitude under computationally economical constraints.

 In the future, we will expand our method to solve singular models, interface problems and anisotropic models. %Time-harmonic Maxwell's equations with singularity}

\end{document}